\documentclass[11pt]{article}
\usepackage{hyperref}
\usepackage{amsmath}
\usepackage{ntheorem}
\usepackage{mathtools}
\usepackage{epsf}
\usepackage{theorem}
\usepackage{indentfirst}
\usepackage{latexsym}
\usepackage{amssymb}
\usepackage[dvips]{graphicx}
\usepackage{flafter}
\usepackage{caption}
\usepackage{textcomp}
\usepackage{supertabular}
\usepackage{longtable, booktabs}
\usepackage{hhline}
\usepackage{bigstrut,bigdelim}
\usepackage{rotating}
\usepackage{multicol}
\usepackage{multirow}
\usepackage{makeidx}
\usepackage{epsfig}
\usepackage{fancyhdr}
\usepackage{threeparttable}
\usepackage{longtable}
\usepackage{mathrsfs}
\usepackage[all]{xy}
\usepackage{xcolor}
\usepackage{cases}
\usepackage{wasysym}
\usepackage{authblk}
\usepackage{authblk}
\usepackage{zhnumber}
\usepackage[bottom]{footmisc}
\newif\ifger
\gerfalse
\newtheorem{definition}{Definition}[section]
\newtheorem{theorem}{Theorem}[section]
\newtheorem{lemma}[definition]{Lemma}
\newtheorem{corollary}[definition]{Corollary}
\newtheorem{remark}{Remark}[section]

\newtheorem*{conjecture}{Conjecture}
\newtheorem{proposition}[definition]{Proposition}
\newtheorem{example}[definition]{Example}
\newtheorem*{funding}{Funding}
\newtheorem*{availability of data and material}{Availability of Data and Material}
\newtheorem*{conflict of interest}{Conflict of interest}

\begin{document}
\baselineskip=19pt

\title{Determination of Some Types of Permutations over $\mathbb{F}_q^2$ \\
with Low-Degree}

\author[a]{Xuan Pang\thanks{E-mail: pangxuan202503@163.com}}
\author[a]{Yangcheng Li\thanks{E-mail: liyc@m.scnu.edu.cn}}
\author[a]{Pingzhi Yuan\thanks{*Corresponding author. E-mail: yuanpz@scnu.edu.cn}}
\author[a]{Yuanpeng Zeng\thanks{E-mail: zengyp2025@163.com}}
\affil[a]{\small \it School of Mathematical Sciences, South China Normal University,\\

 Guangzhou, 510631, P. R. China}

\date{}
\maketitle

\begin{abstract}
The characterization of permutations over finite fields is an important topic in number theory with a long-standing history. This paper presents a systematic investigation of low-degree bivariate polynomial systems $F=(f_1(x,y),f_2(x,y))$ defined over $\mathbb{F}_{q}^2$. Specifically, we employ Hermite's Criterion to completely classify bivariate quadratic permutation polynomial systems, while utilizing the theory of permutation rational functions to give a full classification of bivariate 3-homogeneous permutation polynomial systems.
Furthermore, as an application of our findings, we provide an explicit characterization of the permutation binomials of the form $x^3+ax^{2q+1}$ over $\mathbb{F}_{q^2}$ with characteristic $p\neq3$, thereby resolving a significant special case within this classical research domain.



\end{abstract}

\section{Introduction}
Let $\mathbb{F}_{q^n}$ be the finite field of $q^n$ elements and characteristic $p$, with $\mathbb{F}_{q^n}^*$ be its multiplicative group. A univariate polynomial $f\in\mathbb{F}_{q^n}[x]$ is said to be a {\it permutation polynomial} (PP) if it induces a permutation of $\mathbb{F}_{q^n}$.
Characterizing classes of PPs has been a central and dynamic area of research in finite field theory, driven by their applications in cryptography  \cite{cryp1-RSA,cryp3-Singh}, coding theory  \cite{code1-Ding13,code2-DZ14,code3-Harish}, and  combinatorial designs \cite{cd-DY06}.

Over the past few decades, numerous PPs over  finite fields $\mathbb{F}_{q^n}$ with desirable properties have been found, see \cite{houxd2015FFA,houxd2015book,houxd2023FFA,Wangqiang2025DCC}. However, the basic problem of providing a complete characterization of PPs of a given degree $d$ (or algebraic degree) still remains inherently challenging, primarily due to the sophisticated mathematical computations for verifying permutation properties. To the best of our knowledge, historical progress on PPs' classification is as follows: In 1896, Dickson \cite{Dickson1896} fully classified PPs of degree $d\leq5$ for arbitrary $q$, as well as $d=6$ for odd $q$. Later, in 2010, Li, Chandler and Xiang \cite{Li-degree6-7} extended the classification to degrees $d\in\{6,7\}$ for even $q\geq8$. Subsequent work by Fan \cite{fanxiang2019 degree7,fanxiang2020 degree8odd p, fanxiang2020 degree8even p} provided complete classifications for $d=7$ and $8$. However, no further advancements have been made for higher degrees ($d\geq9$).

 A polynomial over a finite field of the form $f(x)=\sum a_{ij}x^{p^i+p^j}+\sum b_ix^{p^i}+c\in\mathbb{F}_{p^n}[x]$ is classified as a quadratic polynomial because its algebraic degree is at most 2 when interpreted in the additive group structure of $\mathbb{F}_{p^n}$. 
 When $f$ contains only the quadratic terms, that is, $f(x)=\sum a_{ij}x^{p^i+p^j}$, it is specifically termed a Dembowski-Ostrom (DO) polynomial.
 In \cite{Singh-2024-JAA}, Singh and Vishwakarma investigated this class of polynomials and found six quadratic permutations over $\mathbb{F}_{2^n}$, four of which belong to the family of DO polynomials. For more results on DO polynomials, we refer to \cite{Do-poly1,fanxiang2019 degree7}. To date, only few families of such quadratic permutation polynomials are known, which is one of the motivations for our work on classifying quadratic permutations. 

 Notice that the significant challenges in studying permutation polynomials over $\mathbb{F}_{q^n}$, it naturally raises an important question: Could equivalent research be conducted by working directly with $\mathbb{F}_q^n$ instead. This consideration stems from  the vector space isomorphism between $\mathbb{F}_{q^n}$ and $\mathbb{F}_q^n$. As a result, the present study focuses on the investigation of multivariate polynomials. We now introduce the following definition from \cite{Lidl}, which will serve as the foundation for our subsequent analysis.

\begin{definition}\label{def1}
A system of polynomials
$f_1,f_2,\dots,f_n\in\mathbb{F}_{q}[x_1,x_2,\dots,x_n] $
is said to be a permutation of $\mathbb{F}_{q}^n$ if the system of equations
$$\begin{cases}
f_1(x_1,x_2,\dots,x_n)=u_1,\\
~~~~~~\dots\\
f_n(x_1,x_2, \dots,x_n)=u_n,
\end{cases}$$
has exactly one solution in $\mathbb{F}_q^n$ for each $(u_1,u_2,\dots,u_n)\in\mathbb{F}_q^n$.
\end{definition}

In fact, the study of permutation properties over $\mathbb{F}_q^n$ can  be traced back to the research on vectorial Boolean functions (see \cite{boolean1,boolean2}) and APN (Almost Perfect Nonlinear) functions. In \cite{APN1,APN2}, Beierle {\it et al.} constructed two APN permutations in low-dimensional spaces and expressed them in trivariate form $F(x,y,z)=(x^3+uy^2z,y^3+uxz^2,z^3+ux^2y)$. Li and Kaleyski \cite{APN-Likq} later focused on larger dimensions and gave two novel infinite constructions of APN permutations over $\mathbb{F}_{2^m}^3$, with the groundbreaking contribution of representing these infinite families in a generalized trivariate form $F(x,y,z) = (f(x,y,z), g(x,y,z), h(x,y,z))$. Inspired by similar work, Chi, Li and Qu \cite{qu2025} studied  rotatable permutations of $\mathbb{F}_{2^m}^3$ and constructed five infinite classes of 3-homogeneous  rotatable permutations over $\mathbb{F}_{2^m}^3$. While these studies gave some meaningful permutations over the vector space $\mathbb{F}_q^n$, the explicit connection between  $\mathbb{F}_{q^n}$ and $\mathbb{F}_q^n$ has only recently been fully characterized by our work in \cite{algebraic2}, thereby providing essential validation for the equivalence between these two mathematical frameworks. Building upon this foundation, we conclude that studying permutations over $\mathbb{F}_{q}^n$ is indeed a meaningful endeavor.

The rest of this paper is organized as follows. Some notations and preliminaries are stated in Section 2. In Section 3, we first introduce the notions of linear equivalence and  CS equivalence. Based on these preliminaries, we then provide complete characterizations of bivariate quadratic permutations and bivariate 3-homogeneous permutations over finite fields. In particular,  we determine all binomial permutations of the form $x^3+ax^{2q+1}$ over the finite field $\mathbb{F}_{q^2}$ with char$(\mathbb{F}_{q^2})\neq3$, as an application of our results.

\section{Preliminaries}

Let $m$ be a positive integer, a polynomial of the form
$$L(x)=\sum_{i=0}^{m-1}a_ix^{p^i}\in\mathbb{F}_{p^m}[x]$$
is called a linearized polynomials (or a $p$-polynomial) over $\mathbb{F}_{p^m}$. In \cite[Theorem 7.9]{Lidl}, an efficient criterion for determining whether such a polynomial induces a permutation of $\mathbb{F}_{p^n}$, is presented as follows.

\begin{lemma}\label{linear pp}
Let $\mathbb{F}_{q}$ be a finite field of characteristic $p$. Then the linearized polynomial $L(x)$ over $\mathbb{F}_q$ is a permutation polynomial of $\mathbb{F}_q$ if and only if $L(x)$ only has the root zero in $\mathbb{F}_q$.
\end{lemma}

Let $F=(f_1,f_2,\dots, f_n)$ be a system of polynomials with $f_i\in\mathbb{F}_q[x_1,x_2,\dots,x_n]$ for $1\leq i \leq n$. The following lemma provides a necessary and sufficient condition for $F$ to be a  permutation over $\mathbb{F}_q^n$.



\begin{lemma}{\rm{\bf(Hermite's Criterion)}}\label{Hermite}
Let $\mathbb{F}_q$ be of characteristic $p$. Then the system $f_1,\dots,f_n\in \mathbb{F}_{q}[x_1,\dots,x_n]$ induces a permutation of $\mathbb{F}_q^n$ if and only if the following two conditions are satisfied:
\begin{enumerate}
\item[\rm(i)] in the reduction of
$$f_1^{q-1}\dots f_n^{q-1}\bmod(x_1^q-x_1,\dots,x_n^q-x_n)$$
the coefficient of the term $x_1^{q-1}\dots x_n^{q-1}$ is non-zero;
\item[\rm(ii)]in the reduction of
$$f_1^{t_1}\dots f_n^{t_n}\bmod(x_1^q-x_1,\dots,x_n^q-x_n)$$
the coefficient of $x_1^{q-1}\dots x_n^{q-1}$ is zero whenever $t_1,\dots,t_n$ are integers with $0\leq t_i \leq q-1$ for $1\leq i\leq n$, not all $t_i=q-1$, and at least one $t_i\not\equiv0\pmod p$.
\end{enumerate}
\end{lemma}

A fundamental fact regarding permutations over $\mathbb{F}_q^n$ is stated below, which generalizes the univariate case, that is, $f(x)\in\mathbb{F}_q[x]$ permutes $\mathbb{F}_q$ if and only if $f(x)+c$ does so for any $c\in\mathbb{F}_q$.

\begin{lemma}\label{f+contant term}
 Let $f_1,f_2,\dots,f_n\in\mathbb{F}_q[x_1,x_2,\dots,x_n]$, and let $\mathbf{c} = (c_1, \ldots, c_n) \in \mathbb{F}_q^n$ be any constant vector. Then $(f_1,f_2,\dots,f_n)$ induces a permutation of $\mathbb{F}_q^n$ if and only if $(f_1+c_1,f_2+c_2,\dots,f_n+c_n)$ induces a permutation of $\mathbb{F}_q^n$.
 \end{lemma}
\textbf{Proof.} The result follows from the definition of permutation polynomial systems.
 $\hfill\square$

\begin{lemma}\cite{Li&wang}\label{Li&wang}
Let $L(x)$ be a nonzero linearized polynomial over $\mathbb{F}_{2^m}$. Then $x^3+L(x)$ is a permutation polynomial on $\mathbb{F}_{2^m}$ if and only if $m$ is odd, and $L(x)=\theta^2x+\theta x^2$ for some $\theta\in\mathbb{F}_{2^m}^*$.
\end{lemma}


\section{Some classes of permutations over $\mathbb{F}_{q}^2$}
In this section, we give a detailed characterization of several classes of permutations of the vector space $\mathbb{F}_q^2$. Here we first introduce some well-defined terminologies, which will make the statements of our results more concise and explicit.

\begin{definition}
Let $q$ be a prime power and $n$ be a positive integer. Two polynomial systems $F=(f_1,f_2,\dots, f_n)$ and $G=(g_1,g_2,\dots, g_n)$, with $f_i,g_i\in\mathbb{F}_q[x_1,\dots,x_n]$, are called linearly equivalent if there exist nonsingular linear transformations $\rho:\mathbb{F}_q^n\to\mathbb{F}_q^n$ and $\sigma:\mathbb{F}_q^n\to\mathbb{F}_q^n$ such that $G=\rho\circ F\circ \sigma$.
\end{definition}

From this definition, it immediately follows that $F$ is a permutation of $\mathbb{F}_q^n$ if and only if $G$ is a permutation of $\mathbb{F}_q^n$. On the other hand, observe that the identity map $I(x_1,\dots,x_n)=(x_1,\dots,x_n)$ is trivially a permutation system over $\mathbb{F}_{q}^n$. Moreover,  modifying the $i$ component of such polynomial system by adding terms that depend solely on previously determined variables preserves its permutation properties. This motivates the following notion of coordinate shift equivalence:


\begin{definition}
Let $q$ be a prime power and $n$ be a positive integer. Given two polynomial systems $F,G\in(\mathbb{F}_q[x_1,\dots,x_n])^n$ defined as: $$F=(f_1(x_1,\dots,x_k),f_2(x_1,\dots,x_k),\dots, f_k(x_1,\dots,x_k), x_{k+1},\dots, x_n),$$
$$G=(f_1(x_1,\dots,x_k),\dots,f_k(x_1,\dots,x_k),x_{k+1}+h_{k+1}(x_1,\dots,x_k),\dots, x_{n}+h_{n}(x_1,\dots,x_{n-1})),$$
 where $0\leq k\leq n$ and each $h_i\in\mathbb{F}_q[x_1,\dots,x_{i-1}]$. We say $F$ and $G$ are coordinate shift equivalent (abbreviated as CS equivalent). In the special case when $k=0$, the systems reduce to:  $F=(x_1,\dots, x_n)$ and $G=(x_1, x_2+h_2(x_1),\dots, x_n+h_n(x_1,\dots, x_{n-1}))$.
\end{definition}

Clearly, for CS equivalent systems $F$ and $G$, $F$ is a permutation over $\mathbb{F}_q^n$ if and only if $G$ is.


\begin{example}
Let $n=3$, and consider a trivariate polynomial system $G$ over $\mathbb{F}_q^3$ defined by
$$G=(g_1(x_1,x_2,x_3),g_2(x_1,x_2,x_3),g_3(x_1,x_2,x_3))=(x_3,x_1+x_3^2,x_2+x_1x_3).$$
Under the definition of CS equivalence,  $G$ is CS equivalent to $(x_3,x_1,x_2)$.
 If relabeling $(x_3,x_1,x_2)$ as $(y_1,y_2,y_3)$, this system then becomes: $G=(y_1,y_2+y_1^2,y_3+y_1y_2)$, which is CS equivalent to $(y_1,y_2,y_3)$.  Thus, $G$ is CS equivalent to $(x_1, x_2, x_3)$ up to the permutation $\pi = (3\ 1\ 2)$. 
\end{example}

\begin{remark}
The given example shows that $G=(x_3,x_1+x_3^2,x_2+x_1x_3)$ is CS equivalent to both $(x_3,x_1,x_2)$ and $(x_1,x_2,x_3)$ via appropriate variable relabeling. While the original  definition is formulated with a fixed variable ordering, the equivalence naturally extends to permuted coordinates through the transformation $x_i \mapsto x_{\pi(i)}$. Hence, two systems differing only by a permutation of coordinates should be regarded as the same, and this invariance will be implicitly maintained in all subsequent applications of CS equivalence.
\end{remark}

\subsection{quadratic polynomial systems over $\mathbb{F}_q^2$}
 We now initiate our investigation of quadratic polynomial systems. Let us denote $F \sim G$ to signify that two polynomial systems $F$ and $G$ are equivalent under the combined operations of linear equivalences and CS equivalences.

\begin{theorem}\label{th3.1}
Let $\mathbb{F}_q$ be a finite field of characteristic $p$. A quadratic polynomial system $F=(f_1,f_2)$ defined by
\begin{equation}\label{quad-eq}
\begin{cases}
f_1(x,y) = a_1x^2 + a_2xy + a_3y^2 + a_4x + a_5y, \\
f_2(x,y) = b_1x^2 + b_2xy + b_3y^2 + b_4x + b_5y,
\end{cases}
\end{equation}
is a permutation over $\mathbb{F}_q^2$ if and only if:
 \begin{enumerate}
 \item[(i)] for odd $p$, $F\sim(x,y)$,
 \item[(ii)]for even $p$, $F\sim (x,y)$, $(x^2,y)$  $(x,y^2)$, $(x^2,y^2)$ or $(y^2+x, c_1x^2+c_2y^2+c_3x+c_4y)$ with $c_1y^4+(c_2+c_3)y^2+c_4y$ being a linearized permutation over $F_q$.
 \end{enumerate}
\end{theorem}

We prove this theorem in two main parts by separately considering odd and even characteristics. First, we apply Definition \ref{def1} to analyze system (\ref{quad-eq}) and derive necessary and sufficient conditions on the coefficients for permutation behavior over $\mathbb{F}_q^2$. 
Then we further characterize the resulting permutation systems under the two types of equivalence defined above.

Notice that if both $a_2$ and $b_2$ are non-zero, then via a suitable linear transformation, the cross term with coefficient $a_2$ in system (\ref{quad-eq}) can always be eliminated. To streamline our analysis, we therefore assume without loss of generality that $a_2=0$ in system (\ref{quad-eq}).  Henceforth, we will work exclusively with the simplified system where $a_2=0$.

\begin{proposition}\label{pro1}
Let $\mathbb{F}_q$ be a finite field of odd characteristic. The polynomial system $(f_1,f_2)$ defined by
$$
\begin{cases}
f_1(x,y) = a_1x^2 + a_3y^2 + a_4x + a_5y, \\
f_2(x,y) = b_1x^2 + b_2xy + b_3y^2 + b_4x + b_5y,
\end{cases}
$$
is a permutation of $\mathbb{F}_q^2$ if and only if one of the following occurs, up to symmetry in the variables:
\begin{enumerate}
\item[\rm(i)]$f_1=a_5y$ with $a_5\neq0$, $f_2=b_3y^2+b_4x+b_5y$ with $b_4\neq0$;
  \item[\rm(ii)]$f_1=a_3y^2+a_4x+a_5y$ with $a_3\neq0$ and $a_4\neq0$,  $f_2=b_3y^2+b_4x+b_5y$ with $a_3b_4-a_4b_3=0$ and $a_4b_5-a_5b_4\neq0$; 

\item[\rm(iii)] $f_1=a_4x+a_5y$ with $a_4\neq0$ and $a_5\neq0$, $f_2=b_1x^2+b_2xy+b_3y^2+b_4x+b_5y$ with $a_5^2b_1-a_4a_5b_2+b_3a_4^2=b_2a_5-2a_4b_3=0$ and $a_4b_5-a_5b_4\neq0$.
\end{enumerate}
\end{proposition}

\textbf{Proof.}\quad We first deal with the case where both $a_1$ and $a_3$ are non-zero. By completing the square, the polynomial $f_1(x,y)$ can be rewritten as \footnote{The symbol $\triangleq$ denotes ``defined as", meaning the left-hand side is defined as the right-hand side.} 
\begin{equation*}
\begin{aligned}
f_1&=a_1\left(x^2+\frac{a_4}{a_1}x\right)+a_3\left(y^2+\frac{a_5}{a_3}y\right)\\
   &=a_1\left(\left(x+\frac{a_4}{2a_1}\right)^2-\left(\frac{a_4}{2a_1}\right)^2\right)+a_3\left(\left(y+\frac{a_5}{2a_3}\right)^2-\left(\frac{a_5}{2a_3}\right)^2\right)\\
  & \triangleq  a_1X^2+a_3Y^2-\frac{a_4^2}{4a_1}-\frac{a_5^2}{4a_3},
\end{aligned}
\end{equation*}
where the variable substitution is given by $X=x+\frac{a_4}{2a_1}$  and $Y=y+\frac{a_5}{2a_3}$, equivalently, $x=X-\frac{a_4}{2a_1}$ and $y=Y-\frac{a_5}{2a_3}$. Substituting the above variable substitution into $f_2(x,y)$, we obtain
\begin{equation*}
\begin{aligned}
f_2&\triangleq  b_1X^2+b_2XY+b_3Y^2\\
&~~+\left(b_4-\frac{b_1a_4}{a_1}-\frac{b_2a_5}{2a_3}\right)X+\left(b_5-\frac{b_2a_4}{2a_1}-\frac{b_3a_5}{a_3}\right)Y+ a ~constant~term\\
&\triangleq b_1X^2+b_2XY+b_3Y^2+b_{(4)}X+b_{(5)}Y+ a ~constant~term.
\end{aligned}
\end{equation*}
Clearly, $(f_1(x,y),f_2(x,y))$ is a permutation over $\mathbb{F}_{q}^2$ if and only if $(f_1(X,Y),f_2(X,Y))$  is a permutation over $\mathbb{F}_q^2$; by Lemma \ref{f+contant term}, this holds if and only if the system
\begin{equation*}
\begin{cases}
g_1(X,Y)= a_1X^2+a_3Y^2,\\
g_2(X,Y)= b_1X^2+b_2XY+b_3Y^2+b_{(4)}X+b_{(5)}Y,
\end{cases}
\end{equation*}
is a permutation of $\mathbb{F}_q^2$. Let $t_1=q-1,t_2=0$, and  it is easy to verify  that the $X^{q-1}Y^{q-1}$ term in  $g_1^{t_1}g_2^{t_2} \pmod {X^q-X, Y^q-Y}$ is
$$\binom{q-1}{\frac{q-1}{2}}a_1^{\frac{q-1}{2}}a_3^{\frac{q-1}{2}}.$$
By Lucas's theorem \cite[Lemma 2.3]{Lucas-Th}, this coefficient simplifies to $(-a_1a_3)^{\frac{q-1}{2}}\neq0\pmod{p}$ since $a_1,a_3\neq0$. Thus Hermite's criterion dictates that $(g_1,g_2)$ cannot be a permutation of $\mathbb{F}_q^2$. Consequently, if both $a_1$ and $a_3$ are non-zero, the system $(f_1,f_2)$ fails to be a permutation of $\mathbb{F}_q^2$.

We now consider the case where exactly one of $a_1$ and $a_3$ is zero. Suppose $a_1=0$ and $a_3\neq0$ (the case for $a_1\neq 0, a_3=0$ follows by symmetry). By eliminating the $y^2$ term in $f_2$, the system becomes
\begin{equation}\label{eq:10}
\begin{cases}
f_1(x,y) =  a_3y^2 + a_4x + a_5y, \\
f_2(x,y) = b_1x^2 + b_2xy + b_{(4)}x + b_{(5)}y,
\end{cases}
\end{equation}
where $b_{(4)}=b_4-\frac{a_4b_3}{a_3}$, $b_{(5)}=b_5-\frac{a_5b_3}{a_3}$. A simple calculation shows that the coefficient of the $x^{q-1}y^{q-1}$  term in the expansion of $f_1f_2^{q-2} \pmod{x^q-x,y^q-y}$  is
$$a_3 \binom{q-2}{1,q-3,0,0} b_1b_2^{q-3}=(q-2)a_3 b_1b_2^{q-3}.$$
Then Hermite's Criterion implies that at least one of $b_1$ or $b_2$ must vanish. We claim that both $b_1$ and $b_2$ must vanish simultaneously. To see this, suppose for contradiction. Assume $b_1=0$ and $b_2\neq0$, then the coefficient of  $x^{q-1}y^{q-1}$  term in $f_1^0f_2^{q-1}\pmod{x^q-x,y^q-y}$ is $b_2^{q-1}=1\neq0$, violating Hermite's Criterion. Assume $b_1\neq0$ and $b_2=0$, then the coefficient of  $x^{q-1}y^{q-1}$  term in $f_1^{\frac{q-1}{2}}f_2^{\frac{q-1}{2}}\pmod{x^q-x,y^q-y}$ is $a_3^{\frac{q-1}{2}}b_1^{\frac{q-1}{2}}\neq0$ (given $a_3\neq0$), again contradicting Hermite's Criterion. So far we have verified our claim. Therefore, the system (\ref{eq:10}) reduces to
$$
\begin{cases}
f_1(x,y) =  a_3y^2 + a_4x + a_5y, \\
f_2(x,y) =  b_{(4)}x + b_{(5)}y.
\end{cases}
$$
Notably $a_4\neq0$; otherwise $f_1(x,y)$ would degenerate into a quadratic polynomial as $a_3\neq0$, which is incapable of serving as a permutation. By Definition \ref{def1}, it suffices to show that the following equation system
\begin{equation}\label{odd-3}
\begin{cases}
a_3y^2 + a_4x + a_5y=u,\\
b_{(4)}x + b_{(5)}y=v,
\end{cases}
\end{equation}
has exactly one solution in $\mathbb{F}_q^2$ for all $u,v\in\mathbb{F}_q$. Solving $x$ from the first equation in (\ref{odd-3}) gives $x=\frac{u-a_3y^2-a_5y}{a_4}$. Plugging this into the second equation  yields a quadratic equation in $y$:
$$-\frac{a_3b_{(4)}}{a_4}y^2+\left(b_{(5)}-\frac{a_5b_{(4)}}{a_4}\right)y=v-\frac{b_{(4)}}{a_4}u.$$
For unique solvability, this must degenerate to a linear equation, requiring  $b_{(4)}=0$ and $b_{(5)}\neq0$.  Hence  we claim that
$$
\begin{cases}
f_1(x,y) =  a_3y^2 + a_4x + a_5y, \\
f_2(x,y) = b_3y^2+ b_{4}x + b_{5}y.
\end{cases}
$$
with $a_3b_4-a_4b_3=0$ and $a_4b_5-a_5b_4\neq0$ is a permutation of $\mathbb{F}_q^2$. Part (ii) has been proven.

The remaining case is $a_1=a_3=0$, which gives the system:
$$
\begin{cases}
f_1(x,y) = a_4x + a_5y, \\
f_2(x,y) = b_1x^2 + b_2xy + b_3y^2 + b_4x + b_5y.
\end{cases}
$$
It is clear that $a_4$ and $a_5$ cannot vanish simultaneously. For any $u,v\in\mathbb{F}_q$, consider the system $(f_1,f_2)=(u,v)$.  If $a_4=0$ (the case $a_5=0$ is analogous),  part (i) follows immediately. If both $a_4$ and $a_5$ are non-zero, solving for $y$ and substituting gives:
 $$\left(b_1-\frac{b_2a_4}{a_5}+\frac{b_3a_4^2}{a_5^2}\right)x^2+\left(\frac{b_2a_5-2a_4b_3}{a_5^2}u+b_4-\frac{a_4b_5}{a_5}\right)x=v-\frac{b_3}{a_5^2}u^2-\frac{b_5}{a_5}u,$$
from which part (iii) follows.
$\hfill\square$

\medskip

In the following, we focus on the quadratic polynomial systems $F=(f_1,f_2)$ over $\mathbb{F}_{2^m}^2$. A key fact in such fields is that for any $c\in\mathbb{F}_{2^m}$, the equation $ax^2+bx=c$ in $\mathbb{F}_{2^m}$ has a unique solution if and only if either $a=0$ and $b\neq0$ (linear case), or $a\neq0$ and $b=0$ (pure quadratic case). For convenience, we write $P\oplus Q$ to denote the logical exclusive disjunction (XOR), indicating exactly one of $P$ or $Q$ holds.

\begin{proposition}\label{pro2}
Let $\mathbb{F}_q$ be a finite field of even characteristic with $q=2^m$. The polynomial system $(f_1,f_2)$ defined by
$$
\begin{cases}
f_1(x,y) = a_1x^2 + a_3y^2 + a_4x + a_5y, \\
f_2(x,y) = b_1x^2 + b_2xy + b_3y^2 + b_4x + b_5y,
\end{cases}
$$
is a permutation of $\mathbb{F}_q^2$ if and only if one of the following occurs, up to symmetry in the variables:
\begin{enumerate}
\item[\rm(i)]$f_1=a_5y$ with $a_5\neq0$, and either $f_2=b_1x^2+b_3y^2+b_5y$ with $b_1\neq0$ or $f_2=b_3y^2+b_4x+b_5y$ with $b_4\neq0$;
\item[\rm(ii)]$f_1=a_4x+a_5y$ with $a_4\neq0$, $a_5\neq0$, and $f_2=b_1x^2+b_3y^2+b_4x+b_5y$ satisfies $(a_4b_5 + a_5b_4 = 0) \oplus (a_5^2b_1 + a_4^2b_3 = 0)$;
\item[\rm(iii)]$f_1=a_3y^2$ with $a_3\neq0$, and either $f_2=b_1x^2+b_3y^2+b_5y$ with $b_1\neq0$ or $f_2=b_3y^2+b_4x+b_5y$ with $b_4\neq0$;
\item[\rm(iv)]$f_1=a_3y^3+a_4x+a_5y$ with $a_3\neq0,a_4\neq0$, $f_2=b_1x^2+b_3y^2+b_4x+b_5y$ satisfies the following condition:  $L_1(y)=A_1y^4+C_1y^2+D_1y$  has only zero solution in $ \mathbb{F}_{2^m}$, where
 $A_1=\frac{b_1a_3^2}{a_4^2}, C_1=b_3+\frac{b_1a_5^2}{a_4^2}+\frac{b_4a_3}{a_4}$ and  $D_1=\frac{b_4a_5}{a_4}+b_5$;
\item[\rm(v)] $f_1=a_1x^2+a_3y^2$ with $a_1\neq0,a_3\neq0$, and $f_2=b_1x^2+b_3y^2+b_4x+b_5y$ satisfies $(a_1b_3+a_3b_1=0) \oplus (a_1^{1/2}b_5 + a_3^{1/2}b_4=0)$;
\item[\rm(vi)]$f_1= a_1x^2+ a_3y^2 + a_4x + a_5y$  with  $a_1\neq0,a_3\neq0$ and $a_1^{1/2}a_5\neq a_3^{1/2}a_4$, and
         $f_2= b_1x^2+ b_3y^2 + b_4x + b_5y$  satisfies  $L_2(Z)=d_1Z^4+d_2Z^2+d_3Z\in\mathbb{F}_{2^m}[Z]$ having only zero solution in $ \mathbb{F}_{2^m}$, where $d_1=\frac{a_1b_3+a_3b_1}{a_1a_3}$, $d_2=\frac{b_1a_5^2+b_3a_4^2}{a_1a_3}+\left(\frac{a_5}{a_3^{1/2}}+\frac{a_4}{a_1^{1/2}}\right)\left(\frac{b_4}{a_1^{1/2}}+\frac{b_5}{a_3^{1/2}}\right)$, $d_3=\left(\frac{a_4}{a_1^{1/2}}+\frac{a_5}{a_3^{1/2}}\right)\left(\frac{a_4b_5}{(a_1a_3)^{1/2}}+\frac{a_5b_4}{(a_1a_3)^{1/2}}\right)$.
\end{enumerate}
\end{proposition}

\textbf{Proof.}\quad By Definition \ref{def1}, it suffices to show that the quadratic equation system
\begin{numcases}{} 
a_1x^2 + a_3y^2 + a_4x + a_5y = u, \label{eq:1}\\ 
b_1x^2 + b_2xy + b_3y^2 + b_4x + b_5y = v, \label{eq:2} 
\end{numcases}

\noindent has exactly one solution in $\mathbb{F}_{2^m}^2$  for all $u,v\in\mathbb{F}_{2^m}$. We proceed by case analysis based on the coefficients $a_1$ and $a_3$.

 {\bf Case \textup{I}:} Both $a_1$ and $a_3$ are zero. The first equation reduces to a linear equation $a_4x+a_5y=u$, and it is clear that $a_4$ and $a_5$ cannot be zero simultaneously. Hence, we consider two subcases:

Subcase (1.1): Assume that $a_4=0$ and $a_5\neq0$ (the case for $a_4\neq 0$ and $a_5=0$ follows by symmetry), then $y=\frac{u}{a_5}$. Plugging it into Eq. (\ref{eq:2}) and simplifying the equation, we have
\begin{equation}\label{eq:3}
  b_1x^2+\left(\frac{b_2}{a_5}u+b_4\right)x=v+\frac{b_3}{a_5^2}u^2+\frac{b_5 }{a_5}u.
\end{equation}
Now, it is evident that for any $(u,v)\in\mathbb{F}_{2^m}^2$, the quadratic system has a unique solution in $\mathbb{F}_{2^m}^2$ if and only if Eq. ({\ref{eq:3}}) has a unique solution in $\mathbb{F}_{2^m}$. This occurs precisely when Eq. ({\ref{eq:3}}) reduces to either the linear case or the pure quadratic case, that is, $(b_1=0)\oplus (\frac{b_2}{a_5}u+b_4=0)$. Moreover, if $b_1=0$, the condition $\frac{b_2}{a_5}u+b_4\neq0$ for all $u\in\mathbb{F}_{2^m}$ forces $b_2=0$ and $b_4\neq0$; if $b_1\neq0$, the requirement that  $\frac{b_2}{a_5}u+b_4=0$ for all $u$ implies $b_2=b_4=0$. Thus part (i) holds.



Subcase (1.2): Assume that $a_4\neq0$ and $a_5\neq0$. Plugging $x=\frac{u+a_5y}{a_4}$ into Eq. (\ref{eq:2}), we have
\begin{equation}\label{eq:4}
Ay^2+By=v_{(0)}.
\end{equation}
where $A=\frac{b_1a_5^2}{a_4^2}+\frac{b_2a_5}{a_4}+b_3$, $B=\frac{b_2}{a_4}u+\frac{b_4a_5}{a_4}+b_5$, $v_{(0)}=v+\frac{b_1u^2+a_4b_4u}{a_4^2}$.
 It is clear that for any $(u,v)\in\mathbb{F}_q^2$, $(a_4x+a_5y, b_1x^2 + b_2xy + b_3y^2 + b_4x + b_5y)=(u,v)$ has only one solution in $\mathbb{F}_{2^m}^2$ if and only if Eq. (\ref{eq:4}) has one solution in $\mathbb{F}_{2^m}$. Therefore, part (ii) follows from an argument similar to Eq. (\ref{eq:3}).

 If $A=0$, Eq. (\ref{eq:4}) has a unique root $y=\frac{v_{(0)}}{B}$ in $\mathbb{F}_{2^m}$ if and only if $B\neq0$ for all $u\in\mathbb{F}_{2^m}$, which is equivalent to $b_2=0$ and $a_4b_5+a_5b_4\neq0$. Furthermore, when $b_2=0$, $A=0$ yields $a_5^2b_1 + a_4^2b_3 = 0$. If $A\neq0$, Eq. (\ref{eq:4}) has one solution for all $(u,v)\in\mathbb{F}_{2^m}^2$ if and only if $B=0$ for all $u\in\mathbb{F}_{2^m}$. In conclusion, we have proven that when the condition $a_4b_5+a_5b_4\neq0$ and $a_5^2b_1 + a_4^2b_3 = 0$ are satisfied, the system of equations
 \begin{equation}\label{eq:5}
\begin{cases}
 a_4x+a_5y=u,\\
b_1x^2+b_3y^2 + b_4x  + b_5y=v,
\end{cases}
\end{equation}
has only one solution, given by
$$\left(\frac{u}{a_4}+\frac{a_5(a_4^2v+b_1u^2+a_4b_4u)}{a_4^2(a_4b_5+a_5b_4)},\frac{a_4^2v+b_1u^2+a_4b_4u}{a_4(a_4b_5+a_5b_4)}\right).$$

\noindent And the case where $a_4b_5+a_5b_4=0$ and $a_5^2b_1 + a_4^2b_3 \neq 0$, the same system (\ref{eq:5}) has one solution as well. That is,
 $$\left(\frac{u}{a_4}+\frac{a_5}{a_4}\left(\frac{a_4^2v+b_1u^2+a_4b_4u}{b_1 a_5^2+b_3 a_4^2}\right)^{1/2} ,\left(\frac{a_4^2v+b_1u^2+a_4b_4u}{b_1 a_5^2+b_3 a_4^2}\right)^{1/2}\right).$$

As a result, part (ii) follows.
%

 {\bf Case \textup{II}:} Exactly one of $a_1$ and $a_3$ is zero. Assume $a_1=0$ and $a_3\neq0$ (symmetric for $a_1\neq0$ and $a_3=0$). Then Eq. (\ref{eq:1}) reduces to $a_3y^2+a_4x+a_5y=u$.  In the following, we further discuss this based on whether $a_4$ is zero.

Subcase (2.1): For $a_4=0$, we have $a_3y^2+a_5y=u$. Given $a_3\neq0$, this equation has a unique solution if and only if $a_5=0$,  in which case $y=(\frac{u}{a_3})^{1/2}$. Substituting into Eq. (\ref{eq:2}) gives:
$$b_1x^2+\left(\frac{b_2}{a_3^{1/2}}u^{1/2}+b_4\right)x=v+\frac{b_3u+b_5 (a_3u)^{1/2}}{a_3}.$$
To ensure that this equation has exactly one solution for all $(u,v)\in\mathbb{F}_{2^m}$, We use similar reasoning to analyze the possible cases, i.e. $b_1\neq0$ and $b_1=0$. This leads us to conclude that the quadratic system must take one of the forms:
\begin{equation*}
\begin{cases}
 a_3y^2,\\
b_1x^2+b_3y^2+ b_5y ~~{\rm with} ~~b_1\neq0,
\end{cases}
\quad \text{or} \quad
\begin{cases}
 a_3y^2,\\
b_3y^2+b_4x+ b_5y~~{\rm with} ~~b_4\neq0,
\end{cases}
\end{equation*}
corresponding to part (iii) of the theorem.

Subcase (2.2): For $a_4\neq0$,  we solve for $x$ in $f_1(x,y)=u$ to obtain that $x=a_4^{-1}(u+a_3y^2+a_5y)$. Plugging this into Eq. (\ref{eq:2}) gives
\begin{equation}\label{eq:6}
A_1y^4+B_1y^3+C_1y^2+D_1y=v_{1},
\end{equation}
where $A_1=\frac{b_1a_3^2}{a_4^2}$, $B_1=\frac{b_2a_3}{a_4}$, $C_1=b_3+\frac{b_1a_5^2}{a_4^2}+\frac{b_2a_5+b_4a_3}{a_4}$, $D_1=\frac{b_2}{a_4}u+\frac{b_4a_5}{a_4}+b_5$, $v_{1}=v+\frac{b_1u^2+a_4b_4u}{a_4^2}$.
It should be noted that if $B_1\neq0$ (i.e., $b_2\neq0$), then Eq. (\ref{eq:6}) can be rewritten as
$$y^3+\frac{A_1}{B_1}y^4+\frac{C_1}{B_1}y^2+\frac{D_1}{B_1}y=\frac{v_{1}}{B_1}.$$
For this to have exactly one solution for all $(u,v)\in\mathbb{F}_{2^m}^2$, Lemma \ref{Li&wang} implies $A_1=0$ and $C_1^2=B_1D_1$ are necessary. However, since $D_1$ is dependent on the variable $u$, the latter condition cannot hold uniformly for all $u\in\mathbb{F}_{2^m}$. This leads us to consider the alternative case where $B_1=0$. Here, the left-hand side of Eq. (\ref{eq:6}) becomes a linearized polynomial $L_1(y)=A_1y^4+C_1y^2+D_1y$ with $A_1=\frac{b_1a_3^2}{a_4^2}$, $C_1=b_3+\frac{b_1a_5^2}{a_4^2}+\frac{b_4a_3}{a_4}$, $D_1=\frac{b_4a_5}{a_4}+b_5$. According to Lemma \ref{linear pp}, this linearized polynomial has a unique solution for all $(u,v)\in\mathbb{F}_{2^m}^2$ if and only if the only root of it in $\mathbb{F}_{2^m}$ is 0. Hence, part (iv) holds.

{\bf Case \textup{III}: } Both $a_1$ and $a_3$ are non-zero. Owing to the even characteristic of the field, there are $\alpha, \beta\in\mathbb{F}_{2^m}^*$ such that $a_1=\alpha^2$ and $a_3=\beta^2$. Let $X=\alpha x$ and $Y=\beta y$. Then Eq. (\ref{eq:1}) and Eq. ((\ref{eq:2})) can be rewritten as
\begin{equation}\label{eq:9}
\begin{cases}
X^2+Y^2+a_{(4)}X+a_{(5)}Y=u,\\
b_{(1)}X^2+b_{(2)}XY+b_{(3)}Y^2+b_{(4)}X+b_{(5)}Y=v,
\end{cases}
\end{equation}
where $a_{(4)}=a_4/\alpha$,  $a_{(5)}=a_5/\beta$, $b_{(1)}=b_1/\alpha^2$, $b_{(2)}=b_2/\alpha\beta$, $b_{(3)}=b_3/\beta^2$, $b_{(4)}=b_4/\alpha$, $b_{(5)}=b_5/\beta$.

Subcase (3.1): When $a_{(4)}=a_{(5)}$, we have $(X+Y)^2+a_{(4)}(X+Y)=u$. Further, by introducing the substitution
 $Z:=X+Y$, the system (\ref{eq:9}) becomes
\begin{equation*}
\begin{cases}
Z^2+a_{(4)}Z=u,\\
\left(b_{(1)}+b_{(2)}+b_{(3)}\right)Y^2+b_{(2)}YZ+b_{(1)}Z^2+\left(b_{(4)}+b_{(5)}\right)Y+b_{(4)}Z=v.
\end{cases}
\end{equation*}
 The following analysis is analogous to the argument in Subcase (2.1), with details omitted for brevity. This leads to part (v).

Subcase (3.2): When $a_{(4)}\neq a_{(5)}$, we introduce the substitutions
$W:=a_{(4)}X+a_{(5)}Y$ and $Z:=X+Y$, transforming the system (\ref{eq:9}) into:
\begin{equation*}
\begin{cases}
W+Z^2=u,\\
k_1W^2+k_2WZ+k_3Z^2+k_4W+k_5Z=v_{(1)},
\end{cases}
\end{equation*}
where
$$\begin{cases}
k_1=b_{(1)}+b_{(2)}+b_{(3)},\\
k_2=b_{(2)}\left(a_{(4)}+a_{(5)}\right),\\
k_3=b_{(1)}a_{(5)}^2+b_{(2)}a_{(4)}a_{(5)}+b_{(3)}a_{(4)}^2,\\
k_4=\left(a_{(4)}+a_{(5)}\right)\left(b_{(4)}+b_{(5)}\right),\\
k_5=\left(a_{(4)}+a_{(5)}\right)\left(a_{(4)}b_{(5)}+a_{(5)}b_{(4)}\right),\\
v_{(1)}=\left(a_{(4)}^2+a_{(5)}^2\right)v.
\end{cases}$$
Plugging $W=u+Z^2$ into the second equation above, we have
$$k_1Z^4+k_2Z^3+(k_3+k_4)Z^2+(k_2u+k_5)Z=v_{(1)}+k_1u^2+k_4u.$$
A similar discussion with Eq. (\ref{eq:6}) produces the remaining part. This completes the proof.   $\hfill\square$


\begin{remark}
Indeed, the conditions in parts (iv) and (vi) of Proposition \ref{pro2} can be refined through a deeper analysis of cubic equation solutions over finite fields \cite{cubicEq,Tu2018cubic}. For instance, the linearized polynomial $L_1(y)=A_1y^4+C_1y^2+D_1y$ having only the zero solution in $\mathbb{F}_{2^m}$ is equivalent to the cubic equation $A_1y^3+C_1y+D_1=0$ having no non-zero solution in $\mathbb{F}_{2^m}$. While deeper algebraic analysis could further characterize these conditions, we retain the current formulation for conciseness.

\end{remark}

{\bf Proof of Theorem \ref{th3.1}:}
From what we have proven in Propositions \ref{pro1} and \ref{pro2}, we can now easily characterize the obtained permutation systems above in sense of linear equivalence and CS equivalence.

For odd $q$, we have the following:

(i) For $F_1=(a_5y, b_3y^2+b_4x+b_5y)$ with $a_5,b_4\neq0$, it is clear that $F_1$ is CS equivalent to $F_{11}=(a_5y,b_4x)$. Define $\rho_1(x,y)=\left(\frac{y}{b_4},\frac{x}{a_5}\right)$, then we have $\rho_1\circ F_{11}=(x,y)$. Hence $F_1\sim (x,y)$.

 (ii) For $F_2=(a_3y^2+a_4x+a_5y,~b_3y^2+b_4x+b_5y)$ with $a_4\neq0$, $a_3b_4-a_4b_3=0$ and $a_4b_5-a_5b_4\neq0$,  define
 $\rho_2(x,y)=\left(\frac{b_4x-a_4y}{a_5b_4-a_4b_5},\frac{y}{b_4}\right)$ (Note that  if $b_4=0$, then $F_2$ reduces to the case of $F_1$). A direct computation then shows:

  \begin{align*}
    \rho_2 \circ F_2&=\rho_2(a_3y^2+a_4x+a_5y,~b_3y^2+b_4x+b_5y)\\
        &=\left(\frac{b_4(a_3y^2+a_4x+a_5y)-a_4(b_3y^2+b_4x+b_5y)}{a_5b_4-a_4b_5},\frac{b_3y^2+b_4x+b_5y}{b_4}\right)\\
        &= \left(\frac{(a_3b_4-a_4b_3)y^2+(a_5b_4-a_4b_5)y}{a_5b_4-a_4b_5},x+\frac{b_3}{b_4}y^2+\frac{b_5}{b_4}y\right)\\
       &=(y,x+\frac{b_3}{b_4}y^2+\frac{b_5}{b_4}y).
 \end{align*}


(iii) For $F_3=(a_4x+a_5y, b_1x^2+b_2xy+b_3y^2+b_4x+b_5y)$ satisfying $a_5^2b_1-a_4a_5b_2+b_3a_4^2=b_2a_5-2a_4b_3=0$ and $a_4b_5-a_5b_4\neq0$, define $\rho_3(x,y)=\left(\frac{a_5 y}{a_5b_4-a_4b_5},x\right)$ and $\sigma_3(x,y)=\left(x, \frac{-a_4x+y}{a_5}\right)$. Then we can derive that
\begin{align*}
\rho_3 \circ F_3\circ \sigma_3&=\rho_3 \circ F_3\left(x, \frac{-a_4x+y}{a_5}\right)\\
     &=\rho_3 \circ \left(y,~\left(b_1-\frac{a_4b_2}{a_5}+\frac{a_4^2b_3}{a_5^2}\right)x^2+
       \left(\frac{b_2}{a_5}-\frac{2a_4b_3}{a_5^2}\right)xy+\frac{b_3}{a_5^2}y^2+\left(b_4-\frac{a_4b_5}{a_5}\right)x+\frac{b_5}{a_5}y\right)\\
     &=\rho_3 \circ \left(y,~ \frac{b_3}{a_5^2}y^2+\left(b_4-\frac{a_4b_5}{a_5}\right)x+\frac{b_5}{a_5}y\right)\\
    &=\left(x+\frac{b_3}{a_5(a_5b_4-a_4b_5)}y^2+\frac{b_5}{a_5b_4-a_4b_5}y,~y\right).
\end{align*}
It is easily to see that $\left(x+\frac{b_3}{a_5(a_5b_4-a_4b_5)}y^2+\frac{b_5}{a_5b_4-a_4b_5}y,~y\right)$ is CS equivalent to $(x,y)$, and thus $F_3\sim (x,y)$ as desired.

We next handle the even $q$ cases:

(i) Let $G_1=(a_5y,b_1x^2+b_3y^2+b_5y)$ with $a_5,b_1\neq0$, and define the linear transformation $\sigma_1(x,y)=\left(\frac{1}{b_1^{1/2}}y, \frac{1}{a_5}x\right)$. Then we have $G_1\circ \sigma_1=\left(x, y^2+\frac{b_3}{a_5^2}x^2+\frac{b_5}{a_5}x\right)$, which is clearly CS equivalent to $(x,y^2)$. This means $G_1\sim(x,y^2)$. Similarly, it is easily to verify that if $b_4\neq0$, $G_1'=(a_5y, b_3y^2+b_4x+b_5y)$ $\sim (x,y)$.

(ii) Let $G_2=(a_4x+a_5y, b_1x^2+b_3y^2+b_4x+b_5y)$ with $(a_4b_5 + a_5b_4 = 0) \oplus (a_5^2b_1 + a_4^2b_3 = 0)$ and define $\sigma_2(x,y)=\left(\frac{x+a_5y}{a_4}, y\right)$. Then
\begin{align*}
G_2\circ \sigma_2&= \left(x,~  b_1\left(\frac{x}{a_4}+\frac{a_5}{a_4}y \right)^2+b_3y^2+b_4\left(\frac{x}{a_4}+\frac{a_5}{a_4}y\right)+b_5y\right)\\
 &=\left(x,~ \frac{b_1}{a_4^2}x^2+\left(b_3+\frac{b_1a_5^2}{a_4^2}\right)y^2+\frac{b_4}{a_4}x+\left(b_5+\frac{a_5b_4}{a_4}\right)y\right).
\end{align*}
Now, under the conditions where $a_4b_5+a_5b_4=0$ and $a_5^2b_1+a_4^2b_3\neq0$, we immediately derive that $G_1\sim (x,y^2)$. Conversely, when $a_4b_5+a_5b_4\neq0$ and $a_5^2b_1+a_4^2b_3=0$, it follows that $G_1\sim (x,y)$.

(iii) It is clear that  $G_3=(a_3y^2, b_1x^2+b_3y^2+b_5y) \sim (x^2,y^2)$ and $G_3'=(a_3y^2, b_3y^2+b_4x+b_5y)\sim (x^2,y)$, where $b_1\neq0$ in $G_3$ and $b_4\neq0$ in $G_3'$.

(iv) Let $G_4=(a_3y^2+a_4x+a_5y, b_1x^2+b_3y^2+b_4x+b_5y)$ and define $\sigma_4(x,y)=\left(\frac{x}{a_4}+\frac{a_5}{a_4a_3^{1/2}}y, ~ \frac{y}{a_3^{1/2}}\right)$. Then
\begin{align*}
G_4\circ \sigma_4=\left(y^2+x, ~ \frac{b_1}{a_4^2}x^2+\left(\frac{b_1a_5^2}{a_4^2a_3}+\frac{b_3}{a_3}\right)y^2+\frac{b_4}{a_4}x+\left(\frac{a_5b_4}{a_4a_3^{1/2}}+\frac{b_5}{a_3^{1/2}}\right)y\right),
\end{align*}
which induces a permutation over $\mathbb{F}_{2^m}^2$ if and only if the associated  linearized polynomial $L_1'(y)=\frac{b_1}{a_4^2}y^4+\left(\frac{b_1a_5^2}{a_4^2a_3}+\frac{b_3}{a_3}+\frac{b_4}{a_4}\right)y^2+\left(\frac{a_5b_4}{a_4a_3^{1/2}}+\frac{b_5}{a_3^{1/2}}\right)y$ has only the trivial root $y=0$ in $\mathbb{F}_{2^m}$. It is worth noting that this coincides with $L_1(y)$ having only zero solution in $\mathbb{F}_{2^m}$, as seen through the substitution $y_1=a_3^{1/2}y$. Therefore, we have proved that $G_4\sim (y^2+x, c_1x^2+c_2y^2+c_3x+c_4y)$, where $c_1y^4+(c_2+c_3)y^2+c_4y$ is a linear permutation over $\mathbb{F}_{2^m}$.

(v) Let $G_5=(a_1x^2+a_3y^2, b_1x^2+b_3y^2+b_4x+b_5y)$ with $(a_1b_3+a_3b_1=0) \oplus (a_1^{1/2}b_5 + a_3^{1/2}b_4=0)$, and define $\sigma_5(x,y)=\left(\frac{1}{a_1^{1/2}}x+\frac{a_3^{1/2}}{a_1^{1/2}}y,~y\right)$. Then it follows from
$$G_5\circ\sigma_5=\left(x^2,~ \frac{b_1}{a_1}x^2+\left(\frac{a_3b_1}{a_1}+b_3\right)y^2+\frac{b_4}{a_1^{1/2}}x+\left(\frac{a_3^{1/2}b_4}{a_1^{1/2}}+b_5\right)y\right)$$
 that $G_5\sim (x^2,y)$ or $(x^2,y^2)$.

(vi) Let $G_6=(a_1x^2+a_3y^2+a_4x+a_5y, b_1x^2+b_3y^2+b_4x+b_5y)$ satisfy  $L_2(Z)=d_1Z^4+d_2Z^2+d_3Z\in\mathbb{F}_{2^m}[Z]$ having only zero solution in $ \mathbb{F}_{2^m}$. According to the proof above, we easily see that $G_6\sim (y^2+x, c_1x^2+c_2y^2+c_3x+c_4y)$, where $c_1y^4+(c_2+c_3)y^2+c_4y$ is a linear permutation over $\mathbb{F}_{2^m}$.

Thus, the proof of Theorem \ref{th3.1} is hereby completed.
$\hfill\square$

In reviewing the work we have done to fully characterize quadratic permutation systems over $\mathbb{F}_q^2$, we feel the process to be remarkably cumbersome, especially in the case of even characteristic. We believe that this work is truly valuable, which serves as a foundation for subsequent explorations into higher-degrees or multivariate permutations over finite fields. Nevertheless, this is merely a partial breakthrough in the theory of multivariate permutation polynomials. Extending this to dimensions beyond two or degrees higher than quadratic will pose significant challenges, as existing methods struggle with the increased complexity. This is precisely why, in the next section, we only deal with $3$-homogeneous polynomial systems.

\subsection{3-homogeneous polynomial systems over $\mathbb{F}_q^2$}

In this section, we aim to characterize the necessary and sufficient conditions for a bivariate 3-homogeneous polynomial system $F(x,y)=(f_1(x,y),f_2(x,y))$ over $\mathbb{F}_q^2$ to be a permutation. First, note that when $\gcd(3,q-1)\neq1$, $F$ cannot be a permutation of $\mathbb{F}_q^2$  \cite[Proposition 2.7]{qu2025}. Following the methodology employed in Propositions \ref{pro1} and \ref{pro2}, we focus on systems $F$ of the form
\begin{align*}
\begin{cases}
f_1(x,y) = a_1x^3+a_2x^2y + a_3xy^2, \\
f_2(x,y) = b_1x^3+b_2x^2y + b_3xy^2 + b_4y^3.
\end{cases}
\end{align*}
We now claim that $a_1b_4\neq0$ is necessary for $F$ to induce a permutation over $\mathbb{F}_q^2$. To verify this, analyze solutions of $(f_1(x,y,f_2(x,y)))=(0,v)$ in $\mathbb{F}_q^2$ for $v\in\mathbb{F}_q$.  Suppose $a_1=0$. If $b_1b_4\neq0$, then for any $v\in\mathbb{F}_q^*$, the equation system $(f_1,f_2)=(0,v)$ has at least two distinct solutions: $(0,(\frac{v}{b_4})^{1/3})$ and $((\frac{v}{b_1})^{1/3},0)$. If $b_1=0$ (or $b_4=0$), then $(f_1,f_2)=(0,0)$ has at least $q$ distinct solutions of the form $(x,0)$ (or (0,$y$)), where $x\in\mathbb{F}_q$ (or $y\in\mathbb{F}_q$). So $a_1\neq0$. Similarly, suppose $b_4=0$,  and it easy to see that $(x(a_1x^2+a_2xy,a_3y^2), b_1x^3+b_2x^2y + b_3xy^2)=(0,v)$ admits at least $q$ solutions for any $v\in\mathbb{F}_q$, and thus $b_4\neq0$. In conclusion, it suffices to study the case where $F\in(\mathbb{F}_q[x,y])^2$ takes the form below, with $a_1b_4\neq0$ and $\gcd(3,q-1)=1$,
 \begin{align*}
\begin{cases}
f_1(x,y) = a_1x^3+a_2x^2y + a_3xy^2, \\
f_2(x,y) = b_2x^2y + b_3xy^2 + b_4y^3.
\end{cases}
\end{align*}

Before proceeding, we need a few preliminary results. Recall that for a nonconstant rational function $f(X)=\frac{P(X)}{Q(x)}\in\mathbb{F}_q(X)$ with $P(X),Q(X)\in\mathbb{F}_q(X)$  being coprime polynomials, its ${\it degree}$ is defined as deg($f$):=max\{deg($P$),deg($Q$)\}.

\begin{proposition}\label{poly-rational}
Let $f_1,f_2\in\mathbb{F}_q[x,y]$ be coprime homogeneous polynomials of degree $n$, and let  $g(x,y)\in\mathbb{F}_q[x,y]$ be a non-zero homogeneous polynomial of degree $m$. When $m\geq1$, assume further that $(0,0)$ is the only zero of $g(x,y)$ in $\mathbb{F}_q^2$. Suppose these polynomials have the forms:
\begin{align*}
f_1(x,y) &= a_1x^n + a_2x^{n-1}y + \dots+a_nxy^{n-1}, \\
f_2(x,y) &= b_1x^n + b_2x^{n-1}y + \dots+b_nxy^{n-1}+b_{n+1}y^n,\\
g(x,y)&=c_1x^m+c_2x^{m-1}y+\dots+c_mxy^{m-1}+c_{m+1}y^{m},
\end{align*}
If $\gcd(m+n,q-1)=1$, then the following statements are equivalent:

{\rm(i)} the rational function $\frac{f_1(1,t)}{f_2(1,t)}$ is a degree $n$ rational permutation over $\mathbb{P}^1(\mathbb{F}_q):=\mathbb{F}_q\cup \{\infty\}$;

{\rm(ii)} $(f_1(x,y)g(x,y), f_2(x,y)g(x,y))$ induces a permutation over $\mathbb{F}_q^2$.
\end{proposition}

\textbf{Proof.}\quad For the sake of simplicity, let $F(t)=\frac{f_1(1,t)}{f_2(1,t)}=\frac{a_1+a_2t+\dots+a_nt^{n-1}}{b_1+b_2t+\dots+b_nt^{n-1}+b_{n+1}t^n}$, where $t\in\mathbb{P}^1(\mathbb{F}_q)$. When $m\geq1$, we have $g(x,y)\neq0$ for all $(x,y)\neq(0,0)$ since $g(x,y)$ is $m$-homogeneous with $(0,0)$ as its only zero. This implies that both leading coefficients $c_1$ and $c_{m+1}$ are nonzero. When $m = 0$, $g(x,y) = c_0$ is a non-zero constant.

(i)$\Leftarrow$(ii): Suppose that $(f_1(x,y)g(x,y),f_2(x,y)g(x,y))$ induces a permutation over $\mathbb{F}_q^2$. Then

   ~~(a) for any $(u,v)\in\{0\}\times\mathbb{F}_q^*$,
   $(f_1(x,y)g(x,y),f_2(x,y)g(x,y))=(0,v)$ has exactly one solution. This implies $b_{n+1}\neq0$ (as $c_{m+1}\neq0$), ensuring $F(t)$ is a degree $n$ rational function. Consequently, $F(\infty)=0$.

   ~~(b) for any $(u,v)\in\mathbb{F}_q^*\times\{0\}$, $(f_1(x,y)g(x,y),f_2(x,y)g(x,y))=(u,0)$ has exactly one solution. This implies that there exists $(x_u,y_u)\in\mathbb{F}_q^*\times\mathbb{F}_q$ such that
   $$\frac{f_1(x_u,y_u)g(x_u,y_u)}{f_2(x_u,y_u)g(x_u,y_u)}=\frac{x_u^n f_1(1,\frac{y_u}{x_u})}{x_u^n f_2(1,\frac{y_u}{x_u})}=F(t_u)=\infty,$$
where $\frac{u}{0}$ is  interpreted as $\infty$. Thus, $F$ attains $\infty$ at $t_u$.

   ~~(c) for any $(av,v)\in\mathbb{F}_q^*\times\mathbb{F}_q^*$, $(f_1(x,y)g(x,y),f_2(x,y)g(x,y))=(av,v)$ has exactly one solution. This implies that there exists $(x_v,y_v)\in\mathbb{F}_q^*\times\mathbb{F}_q$ such that
     $$\frac{f_1(x_v,y_v)g(x_v,y_v)}{f_2(x_v,y_v)g(x_v,y_v)}=\frac{x_v^n f_1(1,\frac{y_v}{x_v})}{x_v^n f_2(1,\frac{y_v}{x_v})}=F(t_v)=a.$$
Thus for each $a\in\mathbb{F}_q^*$, there exists $t_v\in\mathbb{F}_q\backslash\{t_u\}$ such that $F(t_v)=a$.

Combining (a)-(c), we see that $F(t)$ maps: $\infty$ to $0$, $t_u$ to $\infty$, and covers every $a \in \mathbb{F}_q^*$ at some finite $t_v \neq t_u$. Hence, $F(t):\mathbb{P}^1(\mathbb{F}_q) \to \mathbb{P}^1(\mathbb{F}_q)$ is surjective (and since $\mathbb{P}^1(\mathbb{F}_q)$ is finite, it is also bijective). This indicates that $F(t)$ permutes  $\mathbb{P}^1(\mathbb{F}_q)$.

(i)$\Rightarrow$(ii): Assume that $F(t)$ is a degree $n$ rational permutation over $\mathbb{P}^1(\mathbb{F}_q)$, implying that $b_{n+1}\neq0$. To prove (ii), it suffices to show that  $(f_1(x,y)g(x,y),f_2(x,y)g(x,y))$ runs through all elements of $F_q^2$ as $x$ and $y$ range over $\mathbb{F}_q$. First, consider the case where $(x,y)\in\{0\}\times \mathbb{F}_q$. Substituting $x=0$ into the polynomials, we obtain  $(f_1(0,y)g(0,y),f_2(0,y)g(0,y))=(0,b_{n+1}c_{m+1}y^{m+n})$. Since $\gcd(m+n,q-1)=1$,  the mapping $y\mapsto b_{n+1}c_{m+1}y^{m+n}$ is a permutation of $\mathbb{F}_q$. Thus $(0,b_{n+1}c_{m+1}y^{m+n})$ covers all elements of $\{0\}\times\mathbb{F}_q$ as $y$ ranges over $\mathbb{F}_q$. We next consider the case where $x\neq0$. Clearly, $g(x,y)\neq0$ by hypothesis; let $t=\frac{y}{x}$ ($\neq\infty$). By homogeneity,
$$f_1(x,y)g(x,y)=x^{m+n}f_1(1,t)g(1,t), \quad f_2(x,y)g(x,y)=x^{m+n}f_2(1,t)g(1,t),$$
so $\frac{f_1(x,y)g(x,y)}{f_2(x,y)g(x,y)}=\frac{f_1(1,t)}{f_2(1,t)}=F(t).$ Because $F(t)$ induces a permutation over $\mathbb{P}^1(\mathbb{F}_q)$, there exists some $t_0\in\mathbb{F}_q$ such that $F(t_0)=\infty$, that is, $f_1(1,t_0)\neq0$ and $f_2(1,t_0)=0$. In this case, as $x$ ranges over $\mathbb{F}_q^*$, it is not hard to see that $(x^{m+n}f_1(1,t_0)g(1,t_0),x^{m+n}f_2(1,t_0)g(1,t_0))$ covers $q-1$ elements of $\mathbb{F}_q^*\times\mathbb\{0\}$. For the remaining $t\in\mathbb{F}_q\backslash\{t_0\}$, we have $f_1(1,t)\neq0$ and $f_2(1,t)\neq0$, i.e., $F(t)\neq0$. And in this case, since $\gcd(m+n,q-1)=1$, $$(x_1^{m+n}f_1(1,t_1),x_1^{m+n}f_2(1,t_1))=(x_2^{m+n}f_1(1,t_2),x_2^{m+n}f_2(1,t_2))$$
 would require both $t_1=t_2$ and $x_1=x_2$,  and thus we assert that $(x^{m+n}f_1(1,t),x^{m+n}f_2(1,t))$ runs through $(q-1)^2$ distinct elements $\mathbb{F}_q^*\times\mathbb{F}_q^*$, as $x$ ranges over $\mathbb{F}_q^*$ and $t$ varies over $\mathbb{F}_q\backslash\{t_0\}$. So does $(x^{m+n}f_1(1,t)g(1,t),x^{m+n}f_2(1,t)g(1,t))$ since $g(1,t)\neq0$. In conclusion, $(f_1(x,y)g(x,y),f_2(x,y)g(x,y))$ covers all $q^2$ elements, which implies that $(f_1(x,y)g(x,y),f_2(x,y)g(x,y))$ maps $\mathbb{F}_q$  surjectively onto itself, hence it induces a permutation.
$\hfill\square$

We obtain the following important special case when $m=0$ in Proposition \ref{poly-rational}.
\begin{corollary}\label{cor}
Let $f_1,f_2\in\mathbb{F}_q[x,y]$ be coprime $n$ homogeneous polynomials of the forms:
\begin{align*}
f_1(x,y) &= a_1x^n + a_2x^{n-1}y + \dots+a_nxy^{n-1}, \\
f_2(x,y) &= b_1x^n + b_2x^{n-1}y + \dots+b_nxy^{n-1}+b_{n+1}y^n.
\end{align*}
If $\gcd(n,q-1)=1$, then $(f_1(x,y),f_2(x,y))$ induces a permutation over $\mathbb{F}_q^2$ if and only if $\frac{f_1(1,t)}{f_2(1,t)}$ is a degree $n$ rational permutation over $\mathbb{P}^1(\mathbb{F}_q)$.
\end{corollary}

 On the basis of Proposition \ref{poly-rational} and  Corollary \ref{cor}, we establish a fundamental correspondence between degree $n$  permutation rational functions $f(X)\in\mathbb{F}_q(X)$ and $n$-homogeneous permutation polynomial systems $(f_1(x,y),f_2(x,y))
 \in\mathbb({F}_q[x,y])^2$. This  provides a bridge for simplifying the study of permutation systems over finite fields.

\medskip

Ding and Zieve in \cite[Theorem 1.3]{rational pp} gave a complete characterization of all degree-three rational functions $f(X)\in\mathbb{F}_q(X)$ that permutes $\mathbb{P}^1(\mathbb{F}_q)$.

\begin{lemma}\cite{rational pp}\label{rational pp}
Let $\mu, \nu\in\mathbb{F}_q(X)$ be degree-one  rational functions. A degree-three $f(X)\in\mathbb{F}_q(X)$ permutes $\mathbb{P}^1(\mathbb{F}_q)$ if and only if  one of the following occurs:

{\rm(1)} $\mu\circ X^3 \circ \nu$ where $q\equiv2\pmod3$;

{\rm(2)} $\nu^{-1} \circ X^3 \circ \nu$ where $q\equiv1\pmod3$ and for some $\delta\in\mathbb{F}_{q^2}\backslash\mathbb{F}_{q}$ we have $\nu(X)=(X-\delta^q)/(X-\delta)$ and $\nu^{-1}(X)=(\delta X-\delta^q)$;

{\rm(3)} $\mu\circ (X^3-\gamma X) \circ \nu$ where $3\mid q$ and either $\gamma=0$ or $\gamma$ is a quadratic nonresidue in $\mathbb{F}_q$.
\end{lemma}

We are now ready to state the main result of this section. Recall again that it is enough to examine $3$-homogeneous polynomial systems of the form
$(f_1, f_2) = (a_1x^3 + a_2x^2y + a_3xy^2
, b_2x^2y + b_3xy^2 + b_4y^3),$
where $a_1b_4\neq0$ and $\gcd(3,q-1)=1$.

\begin{theorem}\label{3-hom}
Let $\mathbb{F}_q$ be a finite field with $q\not\equiv1\pmod3$. The system of $3$-homogeneous polynomials defined by
$$
\begin{cases}
f_1(x,y) = a_1x^3 + a_2x^2y + a_3xy^2=xQ_1(x,y),\\
f_2(x,y) = b_2x^2y + b_3xy^2 + b_4y^3=yQ_2(x,y),
\end{cases}
$$
where $a_1b_4\neq0$, $Q_1=a_1x^2+a_2xy+a_3y^2$ and $Q_2=b_2x^2+b_3xy+b_4y^2$, is a permutation of $\mathbb{F}_q^2$  if and only if both $Q_1$ and $Q_2$ are irreducible over $\mathbb{F}_q$, and one of the following occurs:

{\rm(i)} $Q_1(x,y)=k Q_2(x,y)$ for some constant $k\in\mathbb{F}_q^*$;

{\rm(ii)} $\gcd(Q_1(x,y),Q_2(x,y))=1$, and  the rational function $\frac{a_1+a_2t+a_3t^2}{b_2t+b_3t^2+b_4t^3}$ is equivalent to $\mu\circ t^3 \circ \nu$ (for $q\equiv2\pmod 3$) or $\mu\circ (t^3-\gamma t) \circ \nu$ (for $3\mid q$), where $\mu,\nu$ are degree-one rational functions in $\mathbb{F}_q(t)$ and $\gamma\in\mathbb{F}_q$ is either 0 or a  quadratic nonresidue.
\end{theorem}

\textbf{Proof.}\quad  We begin by examining the equation system
\begin{equation}\label{eq:11}
\begin{cases}
 x(a_1x^2 + a_2xy + a_3y^2)=0, \\
 y(b_2x^2 + b_3xy + b_4y^2)=v,
\end{cases}
\end{equation}
to have precisely one solution in $\mathbb{F}_q^2$ for all $v\in\mathbb{F}_q$. Clearly, $(0,(\frac{v}{b_4})^{1/3})$ is always a solution for this system, and thus we must exclude solution with $x\neq0$. When $x\neq0$, set $t=\frac{y}{x}$, so $y=tx$. Substituting into the first equation in (\ref{eq:11}), we have
$x^3( a_1 + a_2t + a_3t^2)=0.$
Since $x\neq0$, we require $a_1 + a_2t + a_3t^2=0$ to have no solutions in $\mathbb{F}_q$. This occurs precisely when the quadratic is irreducible, i,e., $Q_1$ is irreducible.  A similar argument shows that $Q_2$ must be irreducible as well. Therefore, the first part of the theorem has been proven.

If $Q_1(x,y)=k Q_2(x,y)$ for some $k\in\mathbb{F}_q^*$, then by Proposition \ref{poly-rational}, $(f_1(x,y),f_2(x,y))$ permutes $\mathbb{F}_q^2$ since every degree-one ration function permutes $\mathbb{P}^1(\mathbb{F}_q)$ \cite[Lemma 1.2]{rational pp}. Hence, case (i) is trivial. When $\gcd(Q_1(x,y),Q_2(x,y))=1$, observe that $\frac{f_1(1,t)}{f_2(1,t)}=\frac{a_1+a_2t+a_3t^2}{b_2t+b_3t^2+b_4t^3}$, which is a rational function of degree 3 since $b_4\neq0$. Consequently, the remaining part follows from Lemmas \ref{poly-rational} and \ref{rational pp}.
$\hfill\square$

\medskip

We now naturally turn our attention to permutation binomials of the form $f(x)=x^3+ax^{2q+1}$ over $\mathbb{F}_{q^2}$. While the case where $a \in \mathbb{F}_q$ has been completely characterized by Lappano in \cite[Theorem 3.3.2]{Lappano-Ph.D}, we provide an alternative proof using our approaches for completeness. In current paper, we only deal with the case where $q\equiv2\pmod3$. To this end, we first present a clearer characterization of degree 3 permutation rational functions for this particular case as follows.

\begin{proposition}\label{G}
Let $q\equiv2\pmod3$ be a prime power. A degree-three $\frac{g_1}{g_2}=\frac{a_1+a_2t+a_3t^2}{b_1t+b_2t^2+b_3t^3}$ permutes $\mathbb{P}^1(\mathbb{F}_q)$ if and only if it can be expressed as $\frac{1}{d}\cdot\frac{(t - r)^3 - (t - s)^3}{r^3(t - s)^3 - s^3(t - r)^3}$ for some $d,r,s\in\mathbb{F}_q$ with $r\neq s$.
\end{proposition}
\textbf{Proof.}\quad  By Lemma \ref{rational pp}, $\frac{g_1}{g_2}$ permutes $\mathbb{P}^1(\mathbb{F}_q)$ if and only if there exist degree-one $\mu,\nu\in\mathbb{F}_q(t)$ such that $\frac{g_1}{g_2}=\mu\circ t^3\circ\nu$, which is equivalent to $\mu^{-1}\circ\frac{g_1}{g_2}=t^3\circ\nu$. Suppose
 $\mu^{-1}(t)=\frac{p_1t+p_2}{q_1t+q_2}$ and $\nu(t)=\frac{r_1t+r_2}{s_1t+s_2}$, where $p_i,q_i,r_i,s_i\in\mathbb{F}_q$, and $p_1q_2-p_2q_1\neq0$, $r_1s_2-r_2s_1$ are satisfied. Then we have
$$\frac{p_1g_1+p_2g_2}{q_1g_1+q_2g_2}=\frac{(r_1t+r_2)^3}{(s_1t+s_2)^3}.$$
We claim that neither $r_1$ nor $s_1$ is zero. Suppose, for contradiction, that $r_1=0$. Then the right-hand side becomes $\frac{r_2^3}{(s_1t+s_2)^3}$, whose numerator is a constant of degree 0.  By comparing the degrees of the numerators on both sides, this would force $p_1=p_2=0$,  contradicting the fact that $\mu^{-1}$ is degree-one. A similar argument shows $s_1\neq0$. Furthermore, by degree matching again, we easily see that $p_2$ and $q_2$ are non-zero. Hence, the equation can be rewritten as
$$\frac{s_1^3p_2}{r_1^3q_2}\cdot\frac{\frac{p_1}{p_2}g_1+g_2}{\frac{q_1}{q_2}g_1+g_2}=\frac{(t+\frac{r_2}{r_1})^3}{(t+\frac{s_2}{s_1})^3}.$$
By introducing the following relabeled parameters: $$k=\frac{s_1^3p_2}{r_1^3q_2},\quad a=\frac{p_1}{p_2}, \quad  b=\frac{q_1}{q_2}, \quad -r=\frac{r_2}{r_1},\quad -s=\frac{s_2}{s_1},$$ we obtain  a cleaner expression $k\cdot\frac{ag_1+g_2}{bg_1+g_2}=\frac{(t-r)^3}{(t-s)^3}.$
Observing that the leading-order terms (of degree 3) in both the numerator and denominator on the right-hand side are identical, we conclude that $k=1$. A straightforward calculation then yields
$$\frac{g_1}{g_2}=\frac{(t-r)^3-(t-s)^3}{a(t-s)^3-b(t-r)^3}.$$
Finally, notice that $g_2$ has no constant term, we may normalize the coefficients by setting $a=dr^3$ and $b=ds^3$ for some proportionality constant $d\in\mathbb{F}_q^*$. Thus, we arrive at the desired result.
$\hfill\square$

\begin{remark}
By direct calculation,  the discriminants of both the numerator and the denominator in the expression $\frac{1}{d}\cdot\frac{(t - r)^3 - (t - s)^3}{r^3(t - s)^3 - s^3(t - r)^3}$  are equal to $-3$. Notice that $-3$ is a quadratic nonresidue if and only if $q\equiv2\pmod3$. Therefore, for odd prime power $q\equiv2\pmod3$, the irreducibility required in Theorem \ref{3-hom} is automatically satisfied as long as the rational function $\frac{a_1+a_2t+a_3t^2}{b_2t+b_3t^2+b_4t^3}$ takes the form  $\frac{1}{d}\cdot\frac{(t - r)^3 - (t - s)^3}{r^3(t - s)^3 - s^3(t - r)^3}$.
\end{remark}

\begin{proposition}\label{pro3.10}
Let $q$ be an odd prime power with characteristic $p\neq3$, and let $f(x)=x^3+ax^{2q+1}\in\mathbb{F}_{q^2}[x]$ with $a\neq0$. Then $f(x)$ permutes $\mathbb{F}_{q^2}$ if and only if one of the following holds:

(1) If $a\in\mathbb{F}_q$, then

 \quad (1.1) $a=3$ and $q\equiv2\pmod3$;

 \quad (1.2) $a=1$ and $q\equiv 11 \text{ or } 17 \pmod{24}$;

 \quad (1.3) $a=-3$ and $q\equiv -1\pmod{12}$.

(2) If $a\not\in\mathbb{F}_q$, let $u\in\mathbb{F}_q$ be a quadratic nonresidue, and take $\mathbb{F}_{q^2}=\mathbb{F}_q(\alpha)$ where $\alpha^2=u$. Writing $a=a_1+a_2\alpha$ with $a_1,a_2\in\mathbb{F}_q$. Then either:

 \quad (2.1) $a_1=-1$, $a_2=2\sigma/r$ (where $\sigma$ satisfies $z^2-2z-2=0$) and the equation $z^4-4z^2+1=0$ has not roots in $\mathbb{F}_q$; or

 \quad (2.2) 
$a_1 = \frac{3(r^4 + 6ur^2 + u^2)}{(r^2 - u)^2}$ and $a_2 = \frac{ 12(r^2 + u)r}{(r^2 - u)^2}$,
for some arbitrary  $r\in\mathbb{F}_q^*$.
\end{proposition}

\textbf{Proof.}\quad Choose  a quadratic nonresidue $u\in\mathbb{F}_q$, and let $\mathbb{F}_{q^2}=\mathbb{F}_q(\alpha)$, where $\alpha$ satisfies $\alpha^2=u$. For any $x\in\mathbb{F}_{q^2}$, we can express $x$ in the basis $\{1,\alpha\}$ as $x=x_1+x_2\alpha$ with $x_1,x_2\in\mathbb{F}_q$. Similarly, let $a=a_1+a_1\alpha$, with $a_1,a_2\in\mathbb{F}_q$.  Computing the polynomial $f(x)=x^3+ax^{2q+1}$ in this basis yields that
\begin{align*}
f(x)=&(x_1+x_2\alpha)^3+(a_1+a_2\alpha)(x_1+x_2\alpha)^{2q+1}\\
    =&(x_1+x_2\alpha)^3+(a_1+a_2\alpha)(x_1-x_2\alpha)^2(x_1+x_2\alpha)\\
    =&x_1^3+3x_1^2x_2\alpha+3ux_1x_2^2+u\alpha x_2^3+(a_1+a_2\alpha)(x_1^3-ux_1x_2^2-\alpha x_1^2x_2+u\alpha x_2^3)\\
    =&(1+a_1)x_1^3-a_2ux_1^2x_2+(3-a_1)ux_1x_2^2+a_2u^2x_2^3\\
       &+\left(a_2x_1^3+(3-a_1)x_1^2x_2-a_2ux_1x_2^2+(a_1+1)ux_2^3\right)\alpha.
\end{align*}
Now, we define $F=(f_1(x_1,x_2),f_2(x_1,x_2))$ by:
\begin{equation}\label{binomal 3-hom}
\begin{cases}
f_1(x_1,x_2)=(1+a_1)x_1^3-a_2ux_1^2x_2+(3-a_1)ux_1x_2^2+a_2u^2x_2^3,\\
f_2(x_1,x_2)=a_2x_1^3+(3-a_1)x_1^2x_2-a_2ux_1x_2^2+(a_1+1)ux_2^3.
\end{cases}
\end{equation}
It is evident that $f(x)$ permutes $\mathbb{F}_{q^2}$ if and only if $F$ permutes $\mathbb{F}_q^2$. Recall that the condition $(3,q-1)=1$ is necessary for $F$ to be a permutation over $\mathbb{F}_q^2$. Therefore, combined with the fact that char$(\mathbb{F}_q)\neq3$, we only need to consider the case where $q\equiv2\pmod3$ in what follows.

 {\bf Case \textup{1}:} When $a_2=0$ in (\ref{binomal 3-hom}), i.e., $a\in\mathbb{F}_q$, we have
 $$ \begin{cases}
f_1(x_1,x_2)=(1+a_1)x_1^3+(3-a_1)ux_1x_2^2,\\
f_2(x_1,x_2)=(3-a_1)x_1^2x_2+(a_1+1)ux_2^3.
\end{cases}$$

Subcase (1.1): If $a_1=3$, the system reduces to $F=(4x_1^3,4ux_2^3)$. And it is obvious that $F$ permutes $\mathbb{F}_q^2$.

Subcase (1.2): If $a_1=1$, we have $F=(2x_1g(x_1,x_2),2x_2g(x_1,x_2))$, where $g(x_1,x_2)=x_1^2+2ux_2^2$. It is evident that the system $(x_1,x_2)$ permutes $\mathbb{F}_q^2$. Thus by Proposition \ref{poly-rational}, it suffices to examine that $g(x_1,x_2)=0$ has $(0,0)$ as its unique solution. This is equivalent to $-2u$ being a non-square in $\mathbb{F}_q$, which holds precisely when $q\equiv1$ or 3 $\pmod8$. Combined with $q\equiv2\pmod3$, we directly apply the Chinese Remainder Theorem and thus derive $q\equiv11$ or $17\pmod{24}$.

Subcase (1.3): For $a_1\neq1, 3$, let $$G_1(x_1,x_2,a_1,u)=\frac{f_1(x_1,x_2)}{f_2(x_1,x_2)}=\frac{(1+a_1)x_1^3+(3-a_1)ux_1x_2^2}{(3-a_1)x_1^2x_2+(a_1+1)ux_2^3},$$
then we have
$$G=G_1(1,t,a_1,u)=\frac{a_1 + 1 - u(a_1 - 3)t^2}{-(a_1 - 3)t + u(a_1 + 1)t^3},$$
or equivalently
\begin{equation}\label{eq:G}
G=-\frac{(a_1-3)}{(a_1+1)}\cdot\frac{ t^2-\frac{a_1+1}{u(a_1-3)}}{t^3-\frac{a_1-3}{u(a_1+1)}t}.
\end{equation}
By Proposition \ref{G}, the problem now reduces to finding parameters $d,r,s\in\mathbb{F}_q$ such that
$G$ can be expressed as $\frac{1}{d}\cdot\frac{(t - r)^3 - (t - s)^3}{r^3(t - s)^3 - s^3(t - r)^3}$. Expanding this yields the equivalent expression
\begin{equation}\label{expand-coeff} \frac{3\left(t^2-(r+s)t+\frac{r^2+rs+s^2}{3}\right)}{d(r^2+rs+s^2)\left(t^3-\frac{3rs(r+s)}{r^2+rs+s^2}t^2+\frac{3r^2s^2 }{r^2+rs+s^2}t\right)}.
\end{equation}
Comparing this with expression (\ref{eq:G}) leads to the following
$$r+s=0,~~
-\frac{a_1+1}{u(a_1-3)}=\frac{r^2+rs+s^2}{3},$$
$$-\frac{a_1-3}{u(a_1+1)}=\frac{3r^2s^2}{r^2+rs+s^2},~~
-\frac{3rs(r+s)}{r^2+rs+s^2}=0.$$
 A simple computation derives that $\{a_1=-3,r=-s,u=-\frac{1}{s^2}\}$ for arbitrary $s\in\mathbb{F}_q^*$. Since $u$ is a quadratic nonresidue in $\mathbb{F}_q$, $-1$ must be a quadratic nonresidue, implying $q\equiv3\pmod4$. Combining with the earlier condition $q\equiv2\pmod3$ gives $q\equiv{-1}\pmod{12}$.

 {\bf Case \textup{2}:} When $a_2\neq0$ in (\ref{binomal 3-hom}), i.e., $a\not\in\mathbb{F}_q$, we discuss two subcases: $a_1=-1$ and $a_1\neq-1$.

 Subcase (2.1): If $a_1=-1$ in system (\ref{binomal 3-hom}), we have
$$\begin{cases}
f_1(x_1,x_2)=-a_2ux_1^2x_2+4ux_1x_2^2+a_2u^2x_2^3,\\
f_2(x_1,x_2)=a_2x_1^3+4x_1^2x_2-a_2ux_1x_2^2.
\end{cases}.$$
Let
$G_1(x_1,x_2)=\frac{f_1(x_1,x_2)}{f_2(x_1,x_2)}$, and denote $G=G_1(t,1)$, we then obtain
$$G=-\frac{u\left(t^2-\frac{4}{a_2}t-u\right)}{t^3+\frac{4}{a_2}t^2-ut}.$$
For the same reason as Subcase (1.3), we compare its coefficients with those in (\ref{expand-coeff}). This gives the following:
$$r+s=\frac{4}{a_2},~
\frac{r^2+rs+s^2}{3}=-u,~
-\frac{3rs(r+s)}{r^2+rs+s^2}=\frac{4}{a_2},~
\frac{3r^2s^2}{r^2+rs+s^2}=-u.$$
Solving these equations yields  the parametric solution $\{a_2=2\sigma/s, r=(\sigma-3)s, u=(\sigma-3)s^2\}$, where $s\in\mathbb{F}_q^*$ and $\sigma$ satisfies $z^2-2z-2=0$. By assumption that $u$ is a quadratic nonresidue, we have $\sigma-3\in\mathbb{F}_q$ must be a quadratic nonresidue, which indicates that $z^4-4z^2+1=0$ has not roots in $\mathbb{F}_q$.

Subcase (2.1): For $a_1\neq{-1}$, define the transformed system:
\begin{equation*}
\left[
\begin{array}{c}
g_1\\
g_2
\end{array}
\right]
=
\left[
\begin{array}{cc}
a_2 & -(1+a_1)\\
-(1+a_1) & a_2u
\end{array}
\right]
\left[
\begin{array}{c}
f_1\\
f_2
\end{array}
\right].
\end{equation*}
It is evident that the determinant  $a_2^2u-(1+a_1)^2\neq0$, and thus $(f_1,f_2)$ permutes $\mathbb{F}_q^2$ if and only if $(g_1,g_2)$ permutes $\mathbb{F}_q^2$. The explicit forms of $g_1$ and $g_2$ are given by
\begin{equation}\label{eq:14}
\begin{cases}
g_1(x_1,x_2)=-\left(a_2^2u-(a_1-1)^2+4\right)x_1^2x_2+4a_2ux_1x_2^2+\left(a_2^2u-(a_1+1)^2\right)ux_2^3,\\
g_2(x_1,x_2)=\left(a_2^2u-(a_1+1)^2\right)x_1^3+4a_2ux_1^2x_2-u\left(a_2^2u-(a_1-1)^2+4\right)x_1x_2^2.
\end{cases}
\end{equation}

For $-(a_2^2u-(a_1-1)^2+4)=0$, the system (\ref{eq:14}) simplifies to
\begin{equation*}
\begin{cases}
g_1(x_1,x_2)=4a_2ux_1x_2^2+(a_2^2u-(a_1+1)^2)ux_2^3,\\
g_2(x_1,x_2)=(a_2^2u-(a_1+1)^2)x_1^3+4a_2ux_1^2x_2.
\end{cases}
\end{equation*}
We first observe that $a_2^2u-(a_1+1)^2\neq0$ must hold; otherwise, $u$ would be a quadratic residue, contradicting  our assumption. We further claim that this reduced system cannot induce a permutation over $\mathbb{F}_q^2$. To see this, it suffices to examine the solutions of the following system of equation:
\begin{equation*}
\begin{cases}
x_2^2\left(4a_2x_1+(a_2^2u-(a_1+1)^2)x_2\right)=0\\
x_1^2((a_2^2u-(a_1+1)^2)x_1+4a_2ux_2)=v ~~(v\in\mathbb{F}_q^*)
\end{cases}
\end{equation*}
Clearly, $\left(\left(\frac{v}{a_2^2u-(a_1+1)^2}\right)^{1/3},0\right)$ is a solution to this system since $(3,q-1)=1$. On the other hand, it is straightforward to verify that $\left(\frac{(a_1+1)^2-a_2^2u}{4a_2}b,b\right)$ is another solution, where
 $$b^3=\frac{64a_2^3v}{\left(a_2^2u - (a_1+1)^2\right)^2\left(a_2^4u^2 - 2a_1^2a_2^2u + a_1^4 - 4a_1a_2^2u + 4a_1^3 - 18a_2^2u + 6a_1^2 + 4a_1 + 1\right)}\neq0.$$
This shows that the system has at least two distinct solutions, and so cannot be a permutation over $\mathbb{F}_q^2$.

For $ -(a_2^2u-(a_1-1)^2+4)\neq0$, analogously to the previous case, we can derive that
$$G=\frac{-\left(a_2^2u-(a_1-1)^2+4\right) \left(t^{2}-\frac{4 a_2 u}{a_2^2u-(a_1-1)^2+4}t-
\frac{u \left(a_2^{2} u-(a_1+1)^2\right)}{a_2^{2} u-(a_1-1)^2+4}\right)}
{\left(a_2^{2} u-(a_1+1)^2\right) \left(t^{3}+\frac{4 a_2 u }{a_2^{2} u-(a_1+1)^2}t^{2}-
\frac{u \left(a_2^{2} u-(a_1-1)^2+4\right)}{a_2^{2} u-(a_1+1)^2}t\right)}.$$
Further comparison with the expression in (\ref{expand-coeff}) yields $\{a_1=\frac{3(r^2+6rs+s^2)}{(r-s)^2},a_2=\frac{12(r+s)}{(r-s)^2},u=rs\}$ for arbitrary $s\in\mathbb{F}_q^*$.
Substituting $u=rs$ into these expressions, we can rewrite them in terms of $r$ and $u$ as
$$a_1=\frac{3 \left(r^{4}+6 u \,r^{2}+u^{2}\right)}{\left(r^{2}-u\right)^{2}},\quad a_2=\frac{12 \left(u+r^{2}\right) r}{\left(r^{2}-u\right)^{2}}.$$

This completes the proof.
$\hfill\square$

\begin{proposition}\label{pro3.11}
Let $q=2^m$ with odd $m$, and $f(x)=x^3+ax^{2q+1}\in\mathbb{F}_{q^2}[x]$. Then $f(x)$ permutes $\mathbb{F}_{q^2}$ if and only if $a=0$. 
\end{proposition}
\textbf{Proof.}\quad Since $\mathbb{F}_{2^m}$ is a finite field with odd $m$, the polynomial $x^2+x+1$ is irreducible over $\mathbb{F}_{2^m}$. Let $\alpha\in\mathbb{F}_{2^{2m}}$ be a root of $x^2+x+1$, then $\mathbb{F}_{2^{2m}}=\mathbb{F}_{2^m}(\alpha)$. Thus, each $x\in\mathbb{F}_{2^m}$ admits a unique representation $x=x_1+x_2\alpha$ with $x_1,x_2\in\mathbb{F}_{2^m}$, and Similarly, $a=a_1+a_2\alpha$ with $a_1,a_2\in\mathbb{F}_{2^m}$. For brevity, the detailed derivation is included in the appendix, as the rest follows analogously to Proposition \ref{pro3.10}.
$\hfill\square$

From Propositions \ref{pro3.10} and \ref{pro3.11}, we have completely characterized permutation binomials of the form $x^3+ax^{2q+1}$ over $\mathbb{F}_{q^2}$ with characteristic $p\neq3$.

\section{Conclusion and future work}
In this paper, we mainly study bivariate permutations $F=(f_1(x,y),f_2(x,y))$ of low-degree from $\mathbb{F}_{q}^2$ to itself. First, we focus on bivariate quadratic polynomial systems over $\mathbb{F}_q^2$. Using Hermite's Criterion and the definition of multivariate permutations, we obtain the necessary and sufficient conditions on the coefficients for such systems to be permutations over $\mathbb{F}_q^2$. It is worth noting that we introduce two new equivalence relations (linear equivalence and coordinate equivalence) to enhance the clarity of our results. Second, due to the basic correspondence we established between degree $n$ permutation rational functions and bivariate $n$-homogenous permutation polynomial systems, we give a complete characterization for bivariate 3-homogeneous polynomial systems over $\mathbb{F}_q^2$. Finally, by leveraging our results, we determine all  permutation binomials of the form $x^3+ax^{2q+1}\in\mathbb{F}_{q^2}[x]$ with characteristic $p\neq3$.

We believe that our current work  represents a partial breakthrough for multivariate permutations. While bivariate quadratic permutations have been completely classified, it also highlights several fundamental challenges and naturally  gives rise to the following open problem: For $n\geq3$, classify all $n$-variate quadratic polynomial systems over $\mathbb{F}_q^n$. Building on our results, we hereby propose the following conjecture for finite fields of odd characteristic:
\begin{conjecture}
For odd prime power $q$ and $n\geq3$, any quadratic system $F=(f_1,\dots,f_n):\mathbb{F}_{q}^n\to\mathbb{F}_q^n$, where each $f_i\in\mathbb{F}_q[x_1,\dots,x_n]$ is a polynomial of degree at most 2, induces a permutation over $\mathbb{F}_{q}^n$ if and only if $F\sim(x_1,\dots,x_n)$.
\end{conjecture}
This conjecture, if proven, would complete the classification of quadratic permutation polynomials in odd characteristic and provide a foundation for studying more complex systems. However, the even characteristic case, with its distinct challenges, need separate investigation. Additionally, extending the research to degrees higher than quadratic is equally meaningful.

\section*{Acknowledgements}
 \begin{funding}{\rm Pingzhi Yuan was supported by the National Natural Science Foundation of China (Grant No. 12171163) and Guangdong Basic and Applied Basic Research Foundation  (Grant No. 2024A1515010589).
}
\end{funding}

%

\section*{Declarations}
\begin{conflict of interest} {\rm There is no conflict of interest.}
\end{conflict of interest}

\newpage

\section*{Appendix}
\noindent{\bf Proof of Proposition \ref{pro3.11}:}

Let $x=x_1+x_2\alpha$ and $a=a_1+a_2\alpha$ with $x_i,a_i\in\mathbb{F}_{2^m}$, where $\alpha\in\mathbb{F}_{2^{2m}}$ satisfies $\alpha^2+\alpha+1=0$. Then the binomial $f(x)=x^3+ax^{2q+1}$ expands as:
\begin{align*}
f(x)=&(x_1+x_2\alpha)^3+(a_1+a_2\alpha)(x_1+x_2\alpha)^{2q+1}\\
    =&(x_1+x_2\alpha)^3+(a_1+a_2\alpha)(x_1+x_2\alpha)\left(x_1+x_2\alpha^q\right)^{2}\\
    =&(x_1+x_2\alpha)^3+(a_1+a_2\alpha)(x_1+x_2\alpha)\left(x_1+(1+\alpha)x_2\right)^{2}\\
    =&(a_1+1)x_1^3+a_2x_1^2x_2+(a_2+1)x_1x_2^2+(a_1+a_2+1)x_2^3\\
      &+\left(a_2x_1^3+(a_1+a_2+1)x_1^2x_2+(a_1+a_2+_1)x_1x_2^2+a_1x_2^3\right)\alpha.
\end{align*}
Define the polynomials system $F=(f_1(x_1,x_2),f_2(x_1,x_2))$ as follows:
\begin{equation}\label{eq:15}
\begin{cases}
f_1(x_1,x_2)=(a_1+1)x_1^3+a_2x_1^2x_2+(a_2+1)x_1x_2^2+(a_1+a_2+1)x_2^3,\\
f_2(x_1,x_2)=a_2x_1^3+(a_1+a_2+1)x_1^2x_2+(a_1+a_2+1)x_1x_2^2+a_1x_2^3.
\end{cases}
\end{equation}
Now it suffices to consider when $F$ is a permutation of $\mathbb{F}_{2^{m}}^2.$ Notice that $q\equiv2\pmod3$ as $m$ is odd. We proceed by examining two cases: $a_2=0$ and $a_2\neq0$.

 {\bf Case \textup{1}:} When $a_2=0$, i.e., $a\in\mathbb{F}_{2^m}$, the system $F$ reduces to
 $$\begin{cases}
 f_1(x_1,x_2)=(a_1+1)x_1^3+x_1x_2^2+(a_1+1)x_2^3,\\
 f_2(x_1,x_2)=(a_1+1)x_1^2x_2+(a_1+1)x_1x_2^2+a_1x_2^3.
 \end{cases}$$
Using the nondegenerate transformation
$$\begin{bmatrix} g_1 \\ g_2 \end{bmatrix} = \begin{bmatrix} a_1 & -(a_1+1) \\ 0 &  1 \end{bmatrix} \begin{bmatrix} f_1 \\ f_2 \end{bmatrix},
$$
it follows that $(f_1,f_2)$ induces a permutation over $\mathbb{F}_{2^{m}}^2$ if and only if $(g_1,g_2)$ induces a permutation over $\mathbb{F}_{2^{m}}^2$. According to Proposition \ref{poly-rational} and Corollary \ref{cor}, let
$$G(x_1,x_2)=\frac{g_1(x_1,x_2)}{g_2(x_1,x_2)}=
\frac{a_1(a_1+1)x_1^3+(a_1+1)^2x_1^2x_2+(a_1^2+a_1+1)x_1x_2^2}{(a_1+1)x_1^2x_2+(a_1+1)x_1x_2^2+a_1x_2^3},$$
then it follows from $a_1\in\mathbb{F}_{2^m}$ and $a\neq0$ that

$$G(1,t)=\frac{(a_1^2 + a_1 + 1)\cdot\left(t^2+ \frac{(a_1 + 1)^2}{(a_1^2 + a_1 + 1)}t+\frac{a_1(a_1 + 1))}{(a_1^2 + a_1 + 1)}\right)}{a_1\cdot \left(t^3+\frac{a_1 + 1}{a_1}t^2+\frac{a_1 + 1}{a_1}t  \right)}.$$
By Proposition \ref{G}, we now compare coefficients with the expression given in (\ref{expand-coeff}):
$$\frac{(a_1 + 1)^2}{(a_1^2 + a_1 + 1)} = s + r,~~\frac{a_1(a_1 + 1)}{(a_1^2 + a_1 + 1)} = r^2 + rs + s^2,$$
$$\frac{a_1 + 1}{a_1} = \frac{rs(s + r)}{r^2 + rs + s^2},~~~\frac{a_1 + 1}{a_1} = \frac{r^2s^2}{r^2 + rs + s^2}.$$
From the last two equations, we derive the following possible cases:
$$\{a_1 = -1, r = r, s = 0\}, \{a_1 = -1, r = 0, s = s\}, \{a_1 = \frac{s^2 - s + 1}{s - 1}, r =\frac{ s}{s - 1}, s = s\}.$$
By combining these with the first two equations and solving the system, the only possible solution is $\{a_1 = -1, r =s = 0\}$. Indeed, in this case, the system $F=(x_1x_2^2,x_2^3)$, which clearly cannot be a permutation over $\mathbb{F}_{2^m}^2$ since $(x_1x_2^2,x_2^3)=(0,0)$ has at least $q$ distinct solutions in $\mathbb{F}_{2^m}^2$.

 {\bf Case \textup{2}:} When $a_2\neq0$, equivalently, $a\not\in\mathbb{F}_{2^m}$, we consider four distinct subcases.

 Subcase (2.1): If $a_1=1$ in (\ref{eq:15}), the system becomes
 $$\begin{cases}
 f_1(x_1,x_2)=a_2x_1^2x_2+(a_2+1)x_1x_2^2+a_2x_2^3,\\
 f_2(x_1,x_2)=a_2x_1^3+a_2x_1^2x_2+a_2x_1x_2^2+x_2^2.
 \end{cases}
 $$
Let
$$\begin{bmatrix} g_1 \\ g_2 \end{bmatrix} = \begin{bmatrix} 1& a_2 \\ 1 &  0 \end{bmatrix} \begin{bmatrix} f_1 \\ f_2 \end{bmatrix},
$$
and let $G(x_1,x_2)=\frac{g_1(x_1,x_2)}{g_2(x_1,x_2)}$. Then
$$G(1,t)=\frac{(a_2^2+a_2+1)\cdot\left(t^2+\frac{a_2(a_2+1)}{a_2^2+a_2+1}t
  +\frac{a_2^2}{a_2^2+a_2+1}\right)}{a_2\cdot\left(t^3+\frac{a_2+1}{a_2}t^2+t\right)}.$$
By comparing coefficients with (\ref{expand-coeff}), we get $a_2=r=s=1$. However, this solution requires careful verification because $r=s$. Substituting $a_1=a_2=1$ into (\ref{eq:15}) yields
$$\begin{cases}
f_1(x_1,x_2)=x_1^2x_2+x_2^3,\\
f_2(x_1,x_2)=x_1^3+x_1^2x_2+x_1x_2^2+x_2^3.
\end{cases}
$$
It is clear that $(x_1^2x_2+x_2^3,x_1^3+x_1^2x_2+x_1x_2^2+x_2^3)=(0,0)$ holds as long as $x_1=x_2$. In other words, this system has at least $q$ distinct solutions. Thus, this case cannot be a permutation.

Subcase (2.2): If $a_1=0$ in (\ref{eq:15}), the system simplifies to
$$\begin{cases}
f_1(x_1,x_2)=x_1^3+a_2x_1^2x_2+(a_2+1)x_1x_2^2+(a_2+1)x_2^3,\\
f_2(x_1,x_2)=a_2x_1^3+(a_2+1)x_1^2x_2+(a_2+1)x_1x_2^2.
\end{cases}
$$
To further analyze this system, we introduce a nondegenerate linear transformation defined by
$$\begin{bmatrix} g_1 \\ g_2 \end{bmatrix} = \begin{bmatrix} 0& 1 \\ a_2 &  -1 \end{bmatrix} \begin{bmatrix} f_1 \\ f_2 \end{bmatrix},
$$
and similarly, we can obtain
$$G=\frac{(a_2+1)t^2+(a_2+1)t+a_2}{a_2(a_2+1)t^3+(a_2^2+1)t^2+(a_2+a_2+1)t}.$$

For $a_2+1=0$, the system (\ref{eq:15}) takes the form $(f_1,f_2)=(x_1^3+x_1^2x_2, x_1^3)$.
Clearly, it is impossible for this system to induce a permutation of $\mathbb{F}_{2^m}^2$.

 For $a_2+1\neq0$, a similar reasoning to the previous one yields yields $\{a_2=0, r=\beta^2+\beta,s=\beta^2+\beta+1\}$, where $\beta$ is a root of $\beta^2+\beta+1$. However, this contradicts $a\neq0$, and $s\in\mathbb{F}_{2^m}$.

Subcase (2.3): If $a_1+a_2+1=0$ in (\ref{eq:15}), then we have
$$\begin{cases}f_1(x_1,x_2)=(a_1+1)x_1^3+(a_1+1)x_1^2x_2+a_1x_1x_2^2,\\
f_2(x_1,x_2)=(a_1+1)x_1^3+a_1x_2^3,
\end{cases}
$$
and so $$G=\frac{a_1t^2+(a_1+1)t+a_1+1}{a_1t^3+a_1t^2+(a_1+1)t}.$$

For $a_1=0$, we get $a_2=1$, which was already addressed and excluded in Subcase (2.2).

For $a_1\neq0$, via coefficient comparison with (\ref{expand-coeff}), we derive
$$\frac{a_1+1}{a_1}=r+s,~~\frac{a_1+1}{a_1}=r^2+rs+s^2,$$
$$1=\frac{rs(r+s)}{r^2+rs+s^2},~~\frac{a_1+1}{a_1}=\frac{r^2s^2}{r^2+rs+s^2}.$$
Solving $r+s=r^2+rs+s^2=\frac{r^2s^2}{r^2+rs+s^2}$ gives $r=s$, and consequently $a_1=1$, $a_2=0$. This case has already been analyzed in Case 1.

Subcase (2.4): The remaining case is $a_1\neq0$, $a_1+1\neq0$ and $a_1+a_2+1\neq0$. In this case, we carry out a linear transformation for (\ref{eq:14}) and obtain:
$$\begin{bmatrix} g_1 \\ g_2 \end{bmatrix} = \begin{bmatrix} -a_1& 1+a_1+a_2\\ -a_2 &  1+a_1 \end{bmatrix} \begin{bmatrix} f_1 \\ f_2 \end{bmatrix},
$$
that is,
$$\begin{cases}
g_1(x,y)=(a_1^2+a_1a_2+a_2^2+a_1+a_2)x_1^3+(a_1^2+a_1a_2+a_2^2+1)x_1^2 x_2+(a_1^2+a_1a_2+a_2^2+a_1+1)x_1,\\
g_2(x,y)=(a_1^2+a_1a_2+a_2^2+a_2+1)x_1^2x_2+(a_1^2+a_1a_2+a_2^2+1)x_1x_2^2+(a_1^2+a_1a_2+a_2^2+a_1+a_2)x_2^3.
\end{cases}
$$
Note that determinant $-a_1^2+a_1a_2+a_2^2-a_1+a_2\neq0$ is necessary to be a permutation.  We thus proceed by assuming this condition holds. For $a_1^2+a_1a_2+a_2^2+a_1+1\neq0$ and $a_1^2+a_1a_2+a_2^2+a_1+a_2\neq0$, we have
$$G(1,t)=\frac{g_1(1,t)}{g_2(1,t)}=
\frac{(a_1^2+a_1a_2+a_2^2+a_1+1)\left(t^2+\frac{a_1^2+a_1a_2+a_2^2+1}{a_1^2+a_1a_2+a_2^2+a_1+1}t+\frac{a_1^2+a_1a_2+a_2^2+a_1+a_2}{a_1^2+a_1a_2+a_2^2+a_1+1}\right)}
{(a_1^2+a_1a_2+a_2^2+a_1+a_2)\left(t^3+\frac{a_1^2+a_1a_2+a_2^2+1}{a_1^2+a_1a_2+a_2^2+a_1+a_2}t^2+\frac{a_1^2+a_1a_2+a_2^2+a_2+1}{a_1^2+a_1a_2+a_2^2+a_1+a_2}t\right)}.$$
 Comparing the coefficients with expression (\ref{expand-coeff}) yields $a_1,a_2\in\{0,1\}$.  These cases have been previously addressed in our earlier discussion. For $a_1^2+a_1a_2+a_2^2+a_1+1=0$, we easily see that $(g_1,g_2)=(0,v)$ has at least two distinct solutions in $\mathbb{F}_{2^m}^2$, and so cannot be a permutation of $\mathbb{F}_{2^m}^2$. Similarly, for $a_1^2+a_1a_2+a_2^2+a_1+a_2=0$, $(g_1,g_2)=(0,0)$  has at least $q$ solutions in $\mathbb{F}_{2^m}^2$.

  In conclusion, none of these cases can induce a permutation over $\mathbb{F}_{2^m}^2$. This completes the proof of Proposition \ref{pro3.11}.

$\hfill\square$


\end{document}